\documentclass[12pt]{article}
\usepackage{full page, 
amssymb, amscd, graphicx, 
}
    \title{{\bf Logarithmic
tensor category theory, VIII: Braided tensor category
structure on categories of generalized modules for a
conformal vertex algebra}}
    \author{Yi-Zhi Huang, James Lepowsky and Lin Zhang}
    \date{}

\begin{document}
    \bibliographystyle{alpha}
    \maketitle

    \newtheorem{rema}{Remark}[section]
    \newtheorem{propo}[rema]{Proposition}
    \newtheorem{theo}[rema]{Theorem}
   \newtheorem{defi}[rema]{Definition}
    \newtheorem{lemma}[rema]{Lemma}
    \newtheorem{corol}[rema]{Corollary}
     \newtheorem{exam}[rema]{Example}
\newtheorem{assum}[rema]{Assumption}
     \newtheorem{nota}[rema]{Notation}
        \newcommand{\ba}{\begin{array}}
        \newcommand{\ea}{\end{array}}
        \newcommand{\be}{\begin{equation}}
        \newcommand{\ee}{\end{equation}}
        \newcommand{\bea}{\begin{eqnarray}}
        \newcommand{\eea}{\end{eqnarray}}
        \newcommand{\nno}{\nonumber}
        \newcommand{\nn}{\nonumber\\}
        \newcommand{\lbar}{\bigg\vert}
        \newcommand{\p}{\partial}
        \newcommand{\dps}{\displaystyle}
        \newcommand{\bra}{\langle}
        \newcommand{\ket}{\rangle}
 \newcommand{\res}{\mbox{\rm Res}}
\newcommand{\wt}{\mbox{\rm wt}\;}
\newcommand{\swt}{\mbox{\scriptsize\rm wt}\;}
 \newcommand{\pf}{{\it Proof}\hspace{2ex}}
 \newcommand{\epf}{\hspace{2em}$\square$}
 \newcommand{\epfv}{\hspace{1em}$\square$\vspace{1em}}
        \newcommand{\ob}{{\rm ob}\,}
        \renewcommand{\hom}{{\rm Hom}}
\newcommand{\C}{\mathbb{C}}
\newcommand{\R}{\mathbb{R}}
\newcommand{\Z}{\mathbb{Z}}
\newcommand{\N}{\mathbb{N}}
\newcommand{\A}{\mathcal{A}}
\newcommand{\Y}{\mathcal{Y}}
\newcommand{\Arg}{\mbox{\rm Arg}\;}
\newcommand{\comp}{\mathrm{COMP}}
\newcommand{\lgr}{\mathrm{LGR}}

\newcommand{\dlt}[3]{#1 ^{-1}\delta \bigg( \frac{#2 #3 }{#1 }\bigg) }

\newcommand{\dlti}[3]{#1 \delta \bigg( \frac{#2 #3 }{#1 ^{-1}}\bigg) }

 \makeatletter
\newlength{\@pxlwd} \newlength{\@rulewd} \newlength{\@pxlht}
\catcode`.=\active \catcode`B=\active \catcode`:=\active \catcode`|=\active
\def\sprite#1(#2,#3)[#4,#5]{
   \edef\@sprbox{\expandafter\@cdr\string#1\@nil @box}
   \expandafter\newsavebox\csname\@sprbox\endcsname
   \edef#1{\expandafter\usebox\csname\@sprbox\endcsname}
   \expandafter\setbox\csname\@sprbox\endcsname =\hbox\bgroup
   \vbox\bgroup
  \catcode`.=\active\catcode`B=\active\catcode`:=\active\catcode`|=\active
      \@pxlwd=#4 \divide\@pxlwd by #3 \@rulewd=\@pxlwd
      \@pxlht=#5 \divide\@pxlht by #2
      \def .{\hskip \@pxlwd \ignorespaces}
      \def B{\@ifnextchar B{\advance\@rulewd by \@pxlwd}{\vrule
         height \@pxlht width \@rulewd depth 0 pt \@rulewd=\@pxlwd}}
      \def :{\hbox\bgroup\vrule height \@pxlht width 0pt depth
0pt\ignorespaces}
      \def |{\vrule height \@pxlht width 0pt depth 0pt\egroup
         \prevdepth= -1000 pt}
   }
\def\endsprite{\egroup\egroup}
\catcode`.=12 \catcode`B=11 \catcode`:=12 \catcode`|=12\relax
\makeatother

\def\hboxtr{\FormOfHboxtr} 
\sprite{\FormOfHboxtr}(25,25)[0.5 em, 1.2 ex] 

:BBBBBBBBBBBBBBBBBBBBBBBBB |
:BB......................B |
:B.B.....................B |
:B..B....................B |
:B...B...................B |
:B....B..................B |
:B.....B.................B |
:B......B................B |
:B.......B...............B |
:B........B..............B |
:B.........B.............B |
:B..........B............B |
:B...........B...........B |
:B............B..........B |
:B.............B.........B |
:B..............B........B |
:B...............B.......B |
:B................B......B |
:B.................B.....B |
:B..................B....B |
:B...................B...B |
:B....................B..B |
:B.....................B.B |
:B......................BB |
:BBBBBBBBBBBBBBBBBBBBBBBBB |

\endsprite

\def\shboxtr{\FormOfShboxtr} 
\sprite{\FormOfShboxtr}(25,25)[0.3 em, 0.72 ex] 

:BBBBBBBBBBBBBBBBBBBBBBBBB |
:BB......................B |
:B.B.....................B |
:B..B....................B |
:B...B...................B |
:B....B..................B |
:B.....B.................B |
:B......B................B |
:B.......B...............B |
:B........B..............B |
:B.........B.............B |
:B..........B............B |
:B...........B...........B |
:B............B..........B |
:B.............B.........B |
:B..............B........B |
:B...............B.......B |
:B................B......B |
:B.................B.....B |
:B..................B....B |
:B...................B...B |
:B....................B..B |
:B.....................B.B |
:B......................BB |
:BBBBBBBBBBBBBBBBBBBBBBBBB |

\endsprite

\begin{abstract}
This is the eighth part in a series of papers in which we introduce
and develop a natural, general tensor category theory for suitable
module categories for a vertex (operator) algebra.  In this paper
(Part VIII), we construct the braided tensor category structure, using
the previously developed results.
\end{abstract}

\tableofcontents
\vspace{2em}

In this paper, Part VIII of a series of eight papers on logarithmic
tensor category theory, we construct the braided tensor category
structure, using the previously developed results.  The sections,
equations, theorems and so on are numbered globally in the series of
papers rather than within each paper, so that for example equation
(a.b) is the b-th labeled equation in Section a, which is contained in
the paper indicated as follows: In Part I \cite{HLZ1}, which contains
Sections 1 and 2, we give a detailed overview of our theory, state our
main results and introduce the basic objects that we shall study in
this work.  We include a brief discussion of some of the recent
applications of this theory, and also a discussion of some recent
literature.  In Part II \cite{HLZ2}, which contains Section 3, we
develop logarithmic formal calculus and study logarithmic intertwining
operators.  In Part III \cite{HLZ3}, which contains Section 4, we
introduce and study intertwining maps and tensor product bifunctors.
In Part IV \cite{HLZ4}, which contains Sections 5 and 6, we give
constructions of the $P(z)$- and $Q(z)$-tensor product bifunctors
using what we call ``compatibility conditions'' and certain other
conditions.  In Part V \cite{HLZ5}, which contains Sections 7 and 8,
we study products and iterates of intertwining maps and of logarithmic
intertwining operators and we begin the development of our analytic
approach.  In Part VI \cite{HLZ6}, which contains Sections 9 and 10,
we construct the appropriate natural associativity isomorphisms
between triple tensor product functors.  In Part VII \cite{HLZ7},
which contains Section 11, we give sufficient conditions for the
existence of the associativity isomorphisms.  The present paper, Part
VIII, contains Section 12.

\paragraph{Acknowledgments}
The authors gratefully
acknowledge partial support {}from NSF grants DMS-0070800 and
DMS-0401302.  Y.-Z.~H. is also grateful for partial support {}from NSF
grant PHY-0901237 and for the hospitality of Institut des Hautes 
\'{E}tudes Scientifiques in the fall of 2007.

\renewcommand{\theequation}{\thesection.\arabic{equation}}
\renewcommand{\therema}{\thesection.\arabic{rema}}
\setcounter{section}{11}
\setcounter{equation}{0}
\setcounter{rema}{0}

\setcounter{equation}{0}
\setcounter{rema}{0}

\section{The braided tensor category structure}

In this section, we shall complete the formulations and proofs of our
main theorems.  We construct a natural braided monoidal category
structure on the category $\mathcal{C}$. In particular, when
$\mathcal{C}$ is an abelian category, we obtain a natural braided
tensor category structure on $\mathcal{C}$. The strategy and steps in
our construction in this section are essentially the same as those in
\cite{tensorK}, \cite{tensor5} and \cite{rigidity} in the finitely
reductive case but, instead of the corresponding constructions and
results in \cite{tensor1}, \cite{tensor2}, \cite{tensor3} and
\cite{tensor4}, we of course use all the constructions and results we
have obtained in this work except for those in Section 11.  The
present section is independent of Section 11, which provided a method
for verifying the relevant hypotheses.

Our constructions and proofs in this work actually give much more,
namely, the vertex-tensor-categorical structure, in the sense of
\cite{tensorK}, relevant for producing the desired braided tensor
category struture.  We have constructed tensor product bifunctors
depending on a nonzero complex number $z$, along with associativity
isomorphisms between suitable trifunctors constructed from these
bifunctors, and in this section, we shall first give natural
isomorphisms between certain additional functors constructed from
them.  These structures, when enhanced by natural isomorphisms
constructed {}from the Virasoro algebra operators, actually give
vertex tensor category structure (in the sense of \cite{tensorK}).
Our construction of braided tensor category structure in this section
is simply a byproduct of this vertex tensor categorical-type
structure.  In particular, we choose the tensor product bifunctor for
our braided tensor category structure to be the tensor product
bifunctor associated to $z=1$ and we construct all the other necessary
data from the natural isomorphisms mentioned above.  This process of
specializing all our tensor product bifunctors to the bifunctor
associated with only one particular nonzero complex number ``forgets''
all of the essential complex-analytic vertex-tensor-categorical
structure developed in the present work, except for only the
``topological'' information, which is what braided tensor category
structure exhibits.  But even if we are interested only in
constructing our braided tensor category structure, we are still
forced to construct the vertex-tensor-categorical structure first,
because, for instance, iterated tensor products of triples of elements
{\it are not defined in the braided tensor category structure}.

We now return to the setting and assumptions of Section 10, in
particular, Assumption \ref{assum-assoc} and also Remark
\ref{boxtensorchoice}, in which we have made a choice
$\boxtimes_{P(z)}=\hboxtr_{P(z)}'$ of the tensor product bifunctors.
In addition, we also make the following two assumptions, the second of
which implies the convergence condition for intertwining maps in
$\mathcal{C}$ (where we take $W_3 = V$, $M_2=W_4$,
$\mathcal{Y}_{3}=Y_{W_4}$ and $w_{(3)} = {\bf 1}$, and we invoke
Proposition \ref{convergence}; see Definition \ref{conv-conditions}):

\begin{assum}\label{assum-V}
The vertex algebra $V$ as a $V$-module is an object of $\mathcal{C}$.
\end{assum}

\begin{assum}\label{assum-con}
The expansion condition for intertwining maps in $\mathcal{C}$ holds
(see Definition \ref{expansion-conditions}).  Moreover, for objects
$W_1$, $W_2$, $W_3$, $W_4$, $W_{5}$, $M_1$ and $M_{2}$ of
$\mathcal{C}$, logarithmic intertwining operators $\mathcal{Y}_{1}$,
$\mathcal{Y}_{2}$ and $\mathcal{Y}_{3}$ of types ${W_5}\choose
{W_1M_1}$, ${M_1}\choose {W_2M_{2}}$ and ${M_2}\choose {W_3W_{4}}$,
respectively, $z_{1}, z_{2}, z_{3}\in \C^{\times}$ satisfying
$|z_{1}|>|z_{2}|>|z_{3}|>0$, and $w_{(1)}\in W_{1}$, $w_{(2)}\in
W_{2}$, $w_{(3)}\in W_{3}$, $w_{(4)}\in W_{4}$ and $w'_{(5)}\in
W'_{5}$, the series
\begin{equation}\label{3-intw-convp}
\sum_{m, n\in {\mathbb R}}\langle w'_{(5)}, \mathcal{Y}_1(w_{(1)}, z_{1})
\pi_{m}(\mathcal{Y}_2(w_{(2)}, z_{2})\pi_{n}(\mathcal{Y}_3(w_{(3)}, z_{3})
w_{(4)}))\rangle_{W_5}
\end{equation}
is absolutely convergent and can be analytically extended to a
multivalued analytic function on the region given by $z_{1}, z_{2},
z_{3}\ne 0$, $z_{1}\ne z_{2}$, $z_{1}\ne z_{3}$ and $z_{2}\ne z_{3}$,
such that for any set of possible singular points with either
$z_{1}=0$, $z_{2}=0$, $z_{3}=0$, $z_{1}=\infty$, $z_{2}=\infty$,
$z_{3}=\infty$, $z_{1}= z_{2}$, $z_{1}= z_{3}$ or $z_{2}= z_{3}$, this
multivalued analytic function can be expanded near the singularity as
a series having the same form as the expansion near the singular
points of a solution of a system of differential equations with
regular singular points (as defined in Appendix B of \cite{Kn}; recall
Section 11.2).
\end{assum}

\begin{rema}
{\rm By Theorems \ref{sys} and \ref{C_1pp} (see also Remark
\ref{simple-sing}), when $A$ and $\tilde{A}$ are trivial, Assumption
\ref{assum-con} holds if every object of $\mathcal{C}$ satisfies the
$C_{1}$-cofiniteness condition and the quasi-finite dimensionality
condition, or, when $\mathcal{C}$ is in $\mathcal{M}_{sg}$, if every
object of $\mathcal{C}$ is a direct sum of irreducible objects of
$\mathcal{C}$ and there are only finitely many irreducible
$C_{1}$-cofinite objects of $\mathcal{C}$ (up to equivalence).}
\end{rema}

The main result of this work is Theorem \ref{main}, which states that
if $V$ is a M\"{o}bius or conformal vertex algebra and $\mathcal{C}$
is a full subcategory of $\mathcal{M}_{sg}$ or $\mathcal{GM}_{sg}$
satisfying Assumptions \ref{assum-assoc}, \ref{assum-V} and
\ref{assum-con}, then the category $\mathcal{C}$, equipped with the
tensor product bifunctor $\boxtimes$, the unit object $V$, the
braiding isomorphisms $\mathcal{R}$, the associativity isomorphisms
$\mathcal{A}$, and the left and right unit isomorphisms $l$ and $r$,
is an additive braided monoidal category. In particular, if
$\mathcal{C}$ is an abelian category, then equipped with the data
above, it is a braided tensor category.

\subsection{More on tensor products of elements}

Let $W_{1}$, $W_{2}$ and $W_{3}$ be objects of $\mathcal{C}$.  In
Section 7, using the convergence condition, for $z_{1}, z_{2}, z_{3},
z_{4}\in \C^{\times}$ satisfying $|z_{1}|>|z_{2}|>0$ and
$|z_{3}|>|z_{4}|>0$, we have defined the tensor product elements
\[
w_{(1)}\boxtimes_{P(z_{1})}(w_{(2)}\boxtimes_{P(z_{2})}w_{(3)})
\in \overline{W_{1}\boxtimes_{P(z_{1})}(W_{2}\boxtimes_{P(z_{2})}W_{3})}
\]
and
\[
(w_{(1)}\boxtimes_{P(z_{4})}w_{(2)})\boxtimes_{P(z_{3})}w_{(3)}
\in \overline{(W_{1}\boxtimes_{P(z_{4})}W_{2})
\boxtimes_{P(z_{3})}W_{3}},
\]
respectively, for $w_{(1)}\in W_{1}$, $w_{(2)}\in W_{2}$ and
$w_{(3)}\in W_{3}$. In the proof of the commutativity of the hexagon
diagrams below, we shall also need tensor products of elements
$w_{(1)}\in W_{1}$, $w_{(2)}\in W_{2}$ and $w_{(3)}\in W_{3}$ in
$\overline{W_{1}\boxtimes_{P(z_{1})}(W_{2}\boxtimes_{P(z_{2})}W_{3})}$
and in $\overline{(W_{1}\boxtimes_{P(z_{4})}W_{2})
\boxtimes_{P(z_{3})}W_{3}}$ when $z_{1}, z_{2}, z_{3}, z_{4}\in
\C^{\times}$ satisfy $z_{1}\ne z_{2}$ and $z_{3}\ne z_{4}$ but do not
necessarily satisfy the inequality $|z_{1}|>|z_{2}|>0$ or
$|z_{3}|>|z_{4}|>0$. Here we first define these elements.

Let 
$\Y_{1}=\Y_{\boxtimes_{P(z_{1})}, 0}$, 
$\Y_{2}=\Y_{\boxtimes_{P(z_{2})}, 0}$,
$\Y_{3}=\Y_{\boxtimes_{P(z_{3})}, 0}$ and 
$\Y_{4}=\Y_{\boxtimes_{P(z_{4})}, 0}$
be intertwining operators of types 
\[
{W_{1}\boxtimes_{P(z_{1})}(W_{2}\boxtimes_{P(z_{2})}W_{3})\choose 
W_{1}\;\; W_{2}\boxtimes_{P(z_{2})}W_{3}},
\] 
\[
{W_{2}\boxtimes_{P(z_{2})}W_{3}\choose W_{2}\;\; W_{3}},
\]
\[
{(W_{1}\boxtimes_{P(z_{4})}W_{2})\boxtimes_{P(z_{3})}W_{3}\choose 
W_{1}\boxtimes_{P(z_{4})}W_{2}\;\; W_{3}}
\] 
and
\[
{W_{1}\boxtimes_{P(z_{4})}W_{2}\choose W_{1}\;\; W_{2}},
\]
respectively, corresponding to the intertwining maps
$\boxtimes_{P(z_{1})}$, $\boxtimes_{P(z_{2})}$, $\boxtimes_{P(z_{3})}$
and $\boxtimes_{P(z_{4})}$, respectively, as in (\ref{YIp}) and
(\ref{recover}).  Then by Assumption \ref{assum-con},

\[
\langle w', \Y_{1}(w_{(1)}, \zeta_{1})
\Y_{2}(w_{(2)}, \zeta_{2})w_{(3)}\rangle
\]
and 
\[
\langle \tilde{w}', \Y_{3}(\Y_{4}(w_{(1)}, \zeta_{4})
w_{(2)}, \zeta_{3})w_{(3)}\rangle
\]
are absolutely convergent for 
\[
w'\in (W_{1}\boxtimes_{P(z_{1})}(W_{2}\boxtimes_{P(z_{2})}W_{3}))'
\]
and 
\[
\tilde{w}'\in ((W_{1}\boxtimes_{P(z_{4})}W_{2})\boxtimes_{P(z_{3})}W_{3})',
\]
when $|\zeta_{1}|>|\zeta_{2}|>0$ and when $|\zeta_{3}|>|\zeta_{4}|>0$,
respectively, and can be analytically extended to multivalued analytic
functions in the regions given by $\zeta_{1}, \zeta_{2}\ne 0$ and
$\zeta_{1}\ne \zeta_{2}$ and by $\zeta_{3}, \zeta_{4}\ne 0$ and
$\zeta_{3}\ne -\zeta_{4}$, respectively. If we cut these regions along
$\zeta_{1}, \zeta_{2}\ge \R_{+}$ and $\zeta_{3}, \zeta_{4}\in \R_{+}$,
respectively, we obtain simply-connected regions and we can choose
single-valued branches of these multivalued analytic functions.  In
particular, we have the branches of these two multivalued analytic
functions such that their values at points satisfying
$|\zeta_{1}|>|\zeta_{2}|>0$ and $|\zeta_{3}|>|\zeta_{4}|>0$ are
\[
\langle w', \Y_{1}(w_{(1)}, \zeta_{1})
\Y_{2}(w_{(2)}, \zeta_{2})w_{(3)}\rangle
\]
and 
\[
\langle \tilde{w}', \Y_{3}(\Y_{4}(w_{(1)}, \zeta_{4})
w_{(2)}, \zeta_{3})w_{(3)}\rangle,
\]
respectively.

We immediately have:

\begin{propo}
Let $w_{(1)}\in W_{1}$, $w_{(2)}\in W_{2}$ and $w_{(3)}\in W_{3}$.
Then for any $z_{1}, z_{2}, z_{3}, z_{4}\in \C^{\times}$ satisfying
$z_{1}\ne z_{2}$ and $z_{3}\ne -z_{4}$, there exist unique elements
\[
w_{(1)}\boxtimes_{P(z_{1})}(w_{(2)}\boxtimes_{P(z_{2})} w_{(3)}) \in
\overline{W_{1}\boxtimes_{P(z_{1})}(W_{2}\boxtimes_{P(z_{2})}W_{3})}
\]
and 
\[
(w_{(1)}\boxtimes_{P(z_{4})}w_{(2)})\boxtimes_{P(z_{3})} w_{(3)} \in
\overline{(W_{1}\boxtimes_{P(z_{4})}W_{2})\boxtimes_{P(z_{3})}W_{3}}
\]
such that 
for any 
\[
w'\in (W_{1}\boxtimes_{P(z_{1})}(W_{2}\boxtimes_{P(z_{2})}W_{3}))'
\]
and 
\[
\tilde{w}'\in ((W_{1}\boxtimes_{P(z_{4})}W_{2})\boxtimes_{P(z_{3})}W_{3})',
\]
the numbers
\begin{equation}\label{general-tsr-1}
\langle w', w_{(1)}\boxtimes_{P(z_{1})}(w_{(2)}\boxtimes_{P(z_{2})} w_{(3)})
\rangle
\end{equation}
and 
\begin{equation}\label{general-tsr-2}
\langle \tilde{w}', 
(w_{(1)}\boxtimes_{P(z_{4})}w_{(2)})\boxtimes_{P(z_{3})} w_{(3)}
\rangle
\end{equation}
are the values at $(\zeta_{1}, \zeta_{2})=(z_{1}, z_{2})$ and
$(\zeta_{3}, \zeta_{4})=(z_{3}, z_{4})$, respectively, of the branches
of the multivalued analytic functions above of $\zeta_{1}$ and
$\zeta_{2}$ and of $\zeta_{3}$ and $\zeta_{4}$ above,
respectively. \epf
\end{propo}

\begin{rema}
{\rm {}From the definition of 
\[
w_{(1)}\boxtimes_{P(z_{1})}(w_{(2)}\boxtimes_{P(z_{2})} w_{(3)})
\]
and 
\[
(w_{(1)}\boxtimes_{P(z_{4})}w_{(2)})\boxtimes_{P(z_{3})} w_{(3)},
\]
we see that when $|z_{1}|=|z_{2}|$ (with $z_{1}\neq z_{2}$) or $|z_{3}|=|z_{4}|$ 
(with $z_{3}\neq z_{4}$), they are uniquely determined by
\begin{eqnarray*}
\lefteqn{\langle w', w_{(1)}\boxtimes_{P(z_{1})}(w_{(2)}
\boxtimes_{P(z_{2})} w_{(3)})
\rangle}\nn
&&=\lim_{\zeta_{1}\to z_{1},\; \zeta_{2}\to z_{2},\;
|\zeta_{1}|>|\zeta_{2}|>0}
\langle w', \Y_{1}(w_{(1)}, \zeta_{1})\Y_{2}(w_{(2)},
\zeta_{2}) w_{(3)}
\rangle
\end{eqnarray*}
and 
\begin{eqnarray*}
\lefteqn{\langle \tilde{w}', 
(w_{(1)}\boxtimes_{P(z_{4})}w_{(2)})\boxtimes_{P(z_{3})} w_{(3)}
\rangle}\nn
&&=\lim_{\zeta_{3}\to z_{4},\; \zeta_{4}\to z_{4},\;
|\zeta_{3}|>|\zeta_{4}|>0}\langle \tilde{w}', 
\Y_{3}(\Y_{4}(w_{(1)}, \zeta_{4})w_{(2)})\zeta_{3}) w_{(3)}
\rangle
\end{eqnarray*}
for 
\[
w'\in (W_{1}\boxtimes_{P(z_{1})}(W_{2}\boxtimes_{P(z_{2})}W_{3}))'
\]
and 
\[
\tilde{w}'\in ((W_{1}\boxtimes_{P(z_{4})}W_{2})\boxtimes_{P(z_{3})}W_{3})',
\]
where the limits take place in the complex plane with a cut along $\R_{+}$.}
\end{rema}

Propositions \ref{prod=0=>comp=0} and \ref{iter=0=>comp=0},
Corollaries \ref{prospan} and \ref{iterspan} and the definitions of
tensor products of three elements above immediately give:

\begin{propo}
For any $z_{1}, z_{2}, z_{3}, z_{4}\in \C^{\times}$ 
satisfying $z_{1}\ne z_{2}$ and $z_{3}\ne -z_{4}$, the elements of the form 
\begin{eqnarray*}
&\pi_{n}(w_{(1)}\boxtimes_{P(z_{1})}(w_{(2)}\boxtimes_{P(z_{2})}
w_{(3)})),&\\
&\pi_{n}((w_{(1)}\boxtimes_{P(z_{4})}w_{(2)})\boxtimes_{P(z_{3})}
w_{(3)})&
\end{eqnarray*}
for $n\in \R$, $w_{(1)}\in W_{1}$, $w_{(2)}\in W_{2}$, 
$w_{(3)}\in W_{3}$ span 
\begin{eqnarray*}
&W_{1}\boxtimes_{P(z_{1})}(W_{2}\boxtimes_{P(z_{2})}
W_{3}),& \\
&(W_{1}\boxtimes_{P(z_{4})}W_{2})\boxtimes_{P(z_{3})}
W_{3},&
\end{eqnarray*}
respectively.\epf
\end{propo}

Next we discuss tensor products of four elements. These are needed in
the proof of the commutativity of the pentagon diagram below.

\begin{propo}

\begin{enumerate}

\item Let $W_1$, $W_2$, $W_3$, $W_4$, $W_{5}$, $M_1$ and $M_{2}$ be
objects of ${\cal C}$, $z_1,z_2, z_{3}$ nonzero complex numbers
satisfying $|z_1|>|z_2|>|z_{3}|>0$, and $I_1$, $I_2$ and $I_{3}$
$P(z_1)$-, $P(z_2)$- and $P(z_{3})$-intertwining maps of type
${W_5}\choose {W_1M_1}$, ${M_1}\choose {W_2M_{2}}$ and ${M_2}\choose
{W_3W_{4}}$, respectively. Then for $w_{(1)}\in W_1$, $w_{(2)}\in
W_2$, $w_{(3)}\in W_3$, $w_{(4)}\in W_4$ and $w'_{(5)}\in W'_5$, the
series
\begin{equation}\label{3-convp}
\sum_{m, n\in {\mathbb R}}\langle w'_{(5)}, I_1(w_{(1)}\otimes
\pi_m(I_2(w_{(2)}\otimes \pi_n(I_3(w_{(3)}\otimes 
w_{(4)})))))\rangle_{W_5}
\end{equation}
is absolutely convergent. 

\item Let $W_1$, $W_2$, $W_3$, $W_4$, $W_{5}$, $M_3$ and $M_{4}$ be
objects of ${\cal C}$, $z_{1}, z_{23}, z_{3}$ nonzero complex numbers
satisfying $|z_3|>|z_{23}|>0$ and $|z_{1}|>|z_{3}|+|z_{23}|>0$, and $I_1$,
$I_2$ and $I_{3}$ $P(z_{1})$-, $P(z_3)$- and $P(z_{23})$-intertwining
maps of type ${W_5}\choose {W_1M_3}$, ${M_3}\choose {M_{4}W_4}$ and
${M_4}\choose {W_2W_{3}}$, respectively. Then for $w_{(1)}\in W_1$,
$w_{(2)}\in W_2$, $w_{(3)}\in W_3$, $w_{(4)}\in W_4$ and $w'_{(5)}\in
W'_5$, the series
\begin{equation}\label{3-convip}
\sum_{m, n\in {\mathbb R}}\langle w'_{(5)}, I_1(w_{(1)}\otimes
\pi_m(I_2(\pi_n(I_3(w_{(2)}\otimes w_{(3)}))\otimes 
w_{(4)})))\rangle_{W_5}
\end{equation}
is absolutely convergent. 

\item Let $W_1$, $W_2$, $W_3$, $W_4$, $W_{5}$, $M_5$ and $M_{6}$ be
objects of ${\cal C}$, $z_{3}, z_{13}, z_{23}$ nonzero complex numbers
satisfying $|z_{3}|>|z_{13}|>|z_{23}|>0$, and $I_1$, $I_2$ and $I_{3}$
$P(z_{3})$-, $P(z_{13})$- and $P(z_{23})$-intertwining maps of type
${W_5}\choose {M_5W_{4}}$, ${M_5}\choose {W_{1}M_6}$ and ${M_6}\choose
{W_2W_{3}}$, respectively. Then for $w_{(1)}\in W_1$, $w_{(2)}\in
W_2$, $w_{(3)}\in W_3$, $w_{(4)}\in W_4$ and $w'_{(5)}\in W'_5$, the
series
\begin{equation}\label{3-convpi}
\sum_{m, n\in {\mathbb R}}\langle w'_{(5)}, I_1(\pi_m(I_2(w_{(1)}\otimes
\pi_n(I_3(w_{(2)}\otimes w_{(3)}))))\otimes 
w_{(4)})\rangle_{W_5}
\end{equation}
is absolutely convergent. 

\item Let $W_1$, $W_2$, $W_3$, $W_4$, $W_{5}$, $M_7$ and $M_{8}$ be
objects of ${\cal C}$, $z_{3}, z_{23}, z_{12}$ nonzero complex numbers
satisfying $|z_{23}|>|z_{12}|>0$ and $|z_{3}|>|z_{23}|+|z_{12}|>0$, and
$I_1$, $I_2$ and $I_{3}$ $P(z_{3})$-, $P(z_{23})$- and
$P(z_{12})$-intertwining maps of type ${W_5}\choose {M_7W_{4}}$,
${M_7}\choose {M_8W_{3}}$ and ${M_8}\choose {W_1W_{2}}$,
respectively. Then for $w_{(1)}\in W_1$, $w_{(2)}\in W_2$, $w_{(3)}\in
W_3$, $w_{(4)}\in W_4$ and $w'_{(5)}\in W'_5$, the series
\begin{equation}\label{3-convi}
\sum_{m, n\in {\mathbb R}}\langle w'_{(5)}, I_1(\pi_m(I_2(\pi_n(I_3(w_{(1)}\otimes
w_{(2)}))\otimes w_{(3)}))\otimes 
w_{(4)})\rangle_{W_5}
\end{equation}
is absolutely convergent. 

\item Let $W_1$, $W_2$, $W_3$, $W_4$, $W_{5}$, $M_9$ and $M_{10}$ be
objects of ${\cal C}$, $z_{12}, z_2, z_{3}$ nonzero complex numbers
satisfying $|z_2|>|z_{12}|+|z_{3}|>0$, and $I_1$, $I_2$ and $I_{3}$
$P(z_2)$-, $P(z_{12})$- and $P(z_{3})$-intertwining maps of type
${W_5}\choose {M_9M_{10}}$, ${M_9}\choose {W_1W_{2}}$ and
${M_{10}}\choose {W_3W_{4}}$, respectively. Then for $w_{(1)}\in W_1$,
$w_{(2)}\in W_2$, $w_{(3)}\in W_3$, $w_{(4)}\in W_4$ and $w'_{(5)}\in
W'_5$, the series
\begin{equation}\label{3-convcomp}
\sum_{m, n\in {\mathbb R}}\langle w'_{(5)}, I_1(\pi_m(I_2(w_{(1)}\otimes
w_{(2)}))\otimes \pi_n(I_2(w_{(3)}\otimes 
w_{(4)})))\rangle_{W_5}
\end{equation}
is absolutely convergent. 
\end{enumerate}
\end{propo}
\pf The absolute convergence of (\ref{3-convp}) follows immediately
{}from Assumption \ref{assum-con} and Proposition \ref{im:correspond}.

To prove the absolute convergence of (\ref{3-convip}), let $\Y_{1}$,
$\Y_{2}$ and $\Y_{3}$ be the intertwining operators corresponding to
$I_{1}$, $I_{2}$ and $I_{3}$, using $p=0$ as usual, that is,
$\Y_{1}=\Y_{I_{1}, 0}$, $\Y_{2}=\Y_{I_{2}, 0}$ and $\Y_{3}=\Y_{I_{3},
0}$.  We would like to prove that for $w_{(1)}\in W_1$, $w_{(2)}\in
W_2$, $w_{(3)}\in W_3$, $w_{(4)}\in W_4$ and $w'_{(5)}\in W'_5$,
\begin{equation}\label{3-convip-1}
\sum_{m, n\in {\mathbb R}}\langle w'_{(5)}, \Y_1(w_{(1)}, z_{1})
\pi_m(\Y_2(\pi_n(\Y_3(w_{(2)}, z_{23}) w_{(3)}), z_{3}) 
w_{(4)}))\rangle_{W_5}
\end{equation}
is absolutely convergent when $|z_3|>|z_{23}|>0$ and 
$|z_{1}|>|z_{3}|+|z_{23}|>0$. By the $L(0)$-conjugation 
property for intertwining 
operators,  (\ref{3-convip-1}) is equal to 
\begin{eqnarray}\label{3-convip-1.5}
\lefteqn{\sum_{m, n\in {\mathbb R}}\langle e^{(\log z_{3})L'(0)}w'_{(5)}, 
\Y_1(e^{(-\log z_{3})L(0)}w_{(1)}, z_{1}z_{3}^{-1})\cdot}\nn
&&\quad\quad\quad
\cdot\pi_m(\Y_2(\pi_n(\Y_3(e^{-(\log z_{3})L(0)}w_{(2)}, z_{23}z_{3}^{-1}) 
e^{-(\log z_{3})L(0)}w_{(3)}), 1) 
e^{-(\log z_{3})L(0)}w_{(4)})\rangle_{W_5}.\nn
\end{eqnarray}
Since (\ref{3-convip-1}) is equal to (\ref{3-convip-1.5}), we see that
proving that for $w_{(1)}\in W_1$, $w_{(2)}\in W_2$, $w_{(3)}\in
W_3$, $w_{(4)}\in W_4$ and $w'_{(5)}\in W'_5$, (\ref{3-convip-1}) is
absolutely convergent when $|z_3|>|z_{23}|>0$ and
$|z_{1}|>|z_{3}|+|z_{23}|>0$ is equivalent to proving that for
$w_{(1)}\in W_1$, $w_{(2)}\in W_2$, $w_{(3)}\in W_3$, $w_{(4)}\in W_4$
and $w'_{(5)}\in W'_5$,
\begin{equation}\label{3-convip-1.7}
\sum_{m, n\in {\mathbb R}}\langle w'_{(5)}, 
\Y_1(w_{(1)}, \zeta_{1})
\pi_m(\Y_2(\pi_n(\Y_3(w_{(2)}, \zeta_{23}) 
w_{(3)}), 1) 
w_{(4)}))\rangle_{W_5}
\end{equation}
is absolutely convergent when $1>|\zeta_{23}|>0$ and
$|\zeta_{1}|>1+|\zeta_{23}|>0$.

Using Corollary \ref{conv-exp=>asso-op} (the associativity of
intertwining operators), we know that there exist an object $M$ of
${\cal C}$ and intertwining operators $\Y_{4}$ and $\Y_{5}$ of types
${M_{3}\choose W_{2}M}$ and ${M\choose W_{3}W_{4}}$, respectively,
such that for $w'\in M_{3}'$, when $1>|\zeta_{23}|>0$, the series
\[
\sum_{n\in {\mathbb R}}\langle w', 
\Y_2(\pi_n(\Y_3(w_{(2)}, 
\zeta_{23}) w_{(3)}), 1) 
w_{(4)}\rangle_{M_{3}}
=\langle w', 
\Y_2(\Y_3(w_{(2)}, \zeta_{23}) w_{(3)}), 1) 
w_{(4)}\rangle_{M_{3}}
\]
is absolutely convergent and, when $1>|\zeta_{23}|>0$ and
$|\zeta_{23}+1|>1$, its sum is equal to 
\[
\langle w', 
\Y_4(w_{(2)}, \zeta_{23}+1)\Y_{5}(w_{(3)}, 1) 
w_{(4)}\rangle_{M_{3}}.
\]
By Assumption \ref{assum-con}, we know that 
\begin{equation}\label{3-convip-2}
\sum_{m,n\in {\mathbb R}}\langle w'_{(5)}, 
\Y_1(w_{(1)}, \zeta_{1})\pi_{m}(
\Y_4(w_{(2)}, \zeta_{23}+1)
\pi_{n}(\Y_{5}(w_{(3)}, 1) 
w_{(4)}))\rangle_{M_{3}}
\end{equation}
is absolutely convergent and can be analytically extended to a
multivalued analytic function $f(\zeta_{1}, \zeta_{23})$ on the region
given by $\zeta_{1}, \zeta_{23}+1\ne 0$, $\zeta_{1}\ne \zeta_{23}+1$,
$\zeta_{1}\ne 1$ and $\zeta_{23}+1\ne 1$ such that for any set of
possible singular points with either $\zeta_{1}=0$, $\zeta_{23}+1=0$,
$\zeta_{1}= \zeta_{23}+1$, $\zeta_{1}= 1$, $\zeta_{23}=0$,
$\zeta_{1}=\infty$ or $\zeta_{23}+1=\infty$, this multivalued analytic
function can be expanded near its singularity as a series having the
same form as the expansion near the singular points of a solution of a
system of differential equations with regular singular points (again,
as defined in Appendix B of \cite{Kn}; recall Section 11.2).  Since
$(\zeta_{1}, \zeta_{23})=(\infty, 0)$ gives such a set of possible
singular points of $f(\zeta_{1}, \zeta_{23})$, and since for fixed
$\zeta_{1}, \zeta_{23}$ satisfying $|\zeta_{1}|>|\zeta_{23}|+1$ and
$0<|\zeta_{23}|<1$, there exists a positive real number $r<1$ such
that $|\zeta_{1}|>r+1$ and $0<|\zeta_{23}|<r$, $f(\zeta_{1},
\zeta_{23})$ can be expanded as a series in powers of $\zeta_{1}$ and
$\zeta_{23}$ and in nonnegative integral powers of $\log \zeta_{1}$
and $\log \zeta_{23}$ when $|\zeta_{1}|>|\zeta_{23}|+1$ and
$0<|\zeta_{23}|<1$.  Thus we see that (\ref{3-convip-1.7}) is in fact
one value of this expansion of $f(\zeta_{1}, \zeta_{23})$, proving
that (\ref{3-convip-1.7}) is absolutely convergent when
$|\zeta_{1}|>|\zeta_{23}|+1$ and $0<|\zeta_{23}|<1$.

The proofs of the absolute convergence of
(\ref{3-convpi})--(\ref{3-convcomp}) are similar.  \epfv

The special case that the $P(\cdot)$-intertwining maps considered
are $\boxtimes_{P(\cdot)}$ gives the following:

\begin{corol}\label{t-prod-4-elts}
Let $W_1$, $W_2$, $W_3$, $W_4$ be objects of ${\cal C}$ and
$w_{(1)}\in W_1$, $w_{(2)}\in W_2$, $w_{(3)}\in W_3$, $w_{(4)}\in
W_4$.  Then we have:

\begin{enumerate}

\item For $z_1, z_2, z_{3}\in \C^{\times}$ satisfying
$|z_1|>|z_2|>|z_{3}|>0$ and
\[
w'\in (W_{1}\boxtimes_{P(z_{1})} 
(W_{2}\boxtimes_{P(z_{2})} (W_{3}\boxtimes_{P(z_{3})} 
W_{4})))'=W_{1}\hboxtr_{P(z_{1})} 
(W_{2}\boxtimes_{P(z_{2})} (W_{3}\boxtimes_{P(z_{3})} 
W_{4})),
\]
the series 
\begin{equation}\label{4-elts-conv-p}
\sum_{m, n\in {\mathbb R}}\langle w', w_{(1)}\boxtimes_{P(z_{1})} 
\pi_{m}(w_{(2)}\boxtimes_{P(z_{2})} \pi_{n}(w_{(3)}\boxtimes_{P(z_{3})} 
w_{(4)}))\rangle
\end{equation}
is absolutely convergent. 

\item For $z_1, z_{23}, z_3 \in \C^{\times}$  satisfying
$|z_3|>|z_{23}|>0$ and $|z_{1}|>|z_{3}|+|z_{23}|>0$ and
\[
w'\in (W_{1}\boxtimes_{P(z_{1})} 
((W_{2}\boxtimes_{P(z_{23})} W_{3})\boxtimes_{P(z_{3})} 
W_{4}))'=W_{1}\hboxtr_{P(z_{1})} 
((W_{2}\boxtimes_{P(z_{23})} W_{3})\boxtimes_{P(z_{3})} 
W_{4}),
\]
the series 
\begin{equation}\label{4-elts-conv-ip}
\sum_{m, n\in {\mathbb R}}\langle w', w_{(1)}\boxtimes_{P(z_{1})} 
\pi_{m}(\pi_{n}(w_{(2)}\boxtimes_{P(z_{23})} w_{(3)})\boxtimes_{P(z_{3})} 
w_{(4)})\rangle
\end{equation}
is absolutely convergent. 

\item For $z_3, z_{13}, z_{23}\in \C^{\times}$  satisfying
$|z_3|>|z_{13}|>|z_{23}|>0$ and 
\[
w'\in ((W_{1}\boxtimes_{P(z_{13})} 
(W_{2}\boxtimes_{P(z_{23})} W_{3}))\boxtimes_{P(z_{3})} 
W_{4})'=(W_{1}\boxtimes_{P(z_{13})} 
(W_{2}\boxtimes_{P(z_{23})} W_{3}))\hboxtr_{P(z_{3})} 
W_{4},
\]
the series 
\begin{equation}\label{4-elts-conv-pi}
\sum_{m, n\in {\mathbb R}}\langle w', \pi_{m}(w_{(1)}\boxtimes_{P(z_{13})} 
\pi_{n}(w_{(2)}\boxtimes_{P(z_{23})} w_{(3)}))\boxtimes_{P(z_{3})} 
w_{(4)}\rangle
\end{equation}
is absolutely convergent. 

\item For $z_3, z_{23}, z_{12}\in \C^{\times}$  satisfying
$|z_{23}|>|z_{12}|>0$ and $|z_{3}|>|z_{23}|+|z_{12}|>0$ and
\[
w'\in (((W_{1}\boxtimes_{P(z_{12})} 
W_{2})\boxtimes_{P(z_{23})} W_{3})\boxtimes_{P(z_{3})} 
W_{4})'=((W_{1}\boxtimes_{P(z_{12})} 
W_{2})\boxtimes_{P(z_{23})} W_{3})\hboxtr_{P(z_{3})} 
W_{4},
\]
the series 
\begin{equation}\label{4-elts-conv-i}
\sum_{m, n\in {\mathbb R}}\langle w', 
\pi_{m}(\pi_{n}(w_{(1)}\boxtimes_{P(z_{12})} 
w_{(2)})\boxtimes_{P(z_{23})} w_{(3)})\boxtimes_{P(z_{3})} 
w_{(4)}\rangle
\end{equation}
is absolutely convergent. 

\item For $z_{12}, z_2, z_{3}\in \C^{\times}$  satisfying
$|z_2|>|z_{12}|+|z_{3}|>0$ and 
\[
w'\in ((W_{1}\boxtimes_{P(z_{12})} 
W_{2})\boxtimes_{P(z_{2})} (W_{3}\boxtimes_{P(z_{3})} 
W_{4})))'=(W_{1}\boxtimes_{P(z_{12})} 
W_{2})\hboxtr_{P(z_{2})} (W_{3}\boxtimes_{P(z_{3})} 
W_{4})),
\]
the series 
\begin{equation}\label{4-elts-conv-comp}
\sum_{m, n\in {\mathbb R}}\langle w', \pi_{m}(w_{(1)}\boxtimes_{P(z_{12})} 
w_{(2)})\boxtimes_{P(z_{2})} \pi_{n}(w_{(3)}\boxtimes_{P(z_{3})} 
w_{(4)})\rangle
\end{equation}
is absolutely convergent. 
\epf

\end{enumerate}
\end{corol}

Note that the sums of (\ref{4-elts-conv-p}), (\ref{4-elts-conv-ip}),
(\ref{4-elts-conv-pi}), 
(\ref{4-elts-conv-i}), (\ref{4-elts-conv-comp})
define elements of 
\begin{eqnarray}
&\overline{W_{1}\boxtimes_{P(z_{1})}(W_{2}\boxtimes_{P(z_{2})}
(W_{3}\boxtimes_{P(z_{3})}W_{4}))},&\label{tr-4-mod-p}\\
&\overline{W_{1}\boxtimes_{P(z_{1})}((W_{2}\boxtimes_{P(z_{23})}
W_{3})\boxtimes_{P(z_{1})}W_{4})},&\label{tr-4-mod-ip}\\
&\overline{(W_{1}\boxtimes_{P(z_{13})}(W_{2}\boxtimes_{P(z_{23})}
W_{3}))\boxtimes_{P(z_{3})}W_{4}},&\label{tr-4-mod-pi}\\
&\overline{((W_{1}\boxtimes_{P(z_{12})}W_{2})\boxtimes_{P(z_{23})}
W_{3})\boxtimes_{P(z_{3})}W_{4}},&\label{tr-4-mod-i}\\
&\overline{(W_{1}\boxtimes_{P(z_{12})}W_{2})\boxtimes_{P(z_{2})}
(W_{3}\boxtimes_{P(z_{3})}W_{4})},&\label{tr-4-mod-comp}
\end{eqnarray}
respectively, for $z_{1}, z_{2}, z_{3}, z_{12}, z_{13}, z_{23}$
satisfying the corresponding 
inequalities. 
We shall denote these five elements by 
\begin{eqnarray*}
&w_{(1)}\boxtimes_{P(z_{1})}(w_{(2)}\boxtimes_{P(z_{2})}
(w_{(3)}\boxtimes_{P(z_{3})}w_{(4)})), &\\
&w_{(1)}\boxtimes_{P(z_{1})}((w_{(2)}\boxtimes_{P(z_{23})}
w_{(3)})\boxtimes_{P(z_{1})}w_{(4)}),&\\
&(w_{(1)}\boxtimes_{P(z_{13})}(w_{(2)}\boxtimes_{P(z_{23})}
w_{(3)}))\boxtimes_{P(z_{3})}w_{(4)},&\\
&((w_{(1)}\boxtimes_{P(z_{12})}w_{(2)})\boxtimes_{P(z_{23})}
w_{(3)})\boxtimes_{P(z_{3})}w_{(4)},&\\
&(w_{(1)}\boxtimes_{P(z_{12})}w_{(2)})\boxtimes_{P(z_{2})}
(w_{(3)}\boxtimes_{P(z_{3})}w_{(4)}),&
\end{eqnarray*}
respectively.

\begin{propo}
The elements of the form 
\begin{eqnarray*}
&\pi_{n}(w_{(1)}\boxtimes_{P(z_{1})}(w_{(2)}\boxtimes_{P(z_{2})}
(w_{(3)}\boxtimes_{P(z_{3})}w_{(4)}))),&\\
&\pi_{n}(w_{(1)}\boxtimes_{P(z_{1})}((w_{(2)}\boxtimes_{P(z_{23})}
w_{(3)})\boxtimes_{P(z_{1})}w_{(4)})),&\\
&\pi_{n}((w_{(1)}\boxtimes_{P(z_{13})}(w_{(2)}\boxtimes_{P(z_{23})}
w_{(3)}))\boxtimes_{P(z_{3})}w_{(4)}),&\\
&\pi_{n}(((w_{(1)}\boxtimes_{P(z_{12})}w_{(2)})\boxtimes_{P(z_{23})}
w_{(3)})\boxtimes_{P(z_{3})}w_{(4)}),&\\
&\pi_{n}((w_{(1)}\boxtimes_{P(z_{12})}w_{(2)})\boxtimes_{P(z_{2})}
(w_{(3)}\boxtimes_{P(z_{3})}w_{(4)}))&
\end{eqnarray*}
for $n\in \R$, $w_{(1)}\in W_{1}$, $w_{(2)}\in W_{2}$, 
$w_{(3)}\in W_{3}$, $w_{(4)}\in W_{4}$ span 
\begin{eqnarray*}
&W_{1}\boxtimes_{P(z_{1})}(W_{2}\boxtimes_{P(z_{2})}
(W_{3}\boxtimes_{P(z_{3})}W_{4})),&\\
&W_{1}\boxtimes_{P(z_{1})}((W_{2}\boxtimes_{P(z_{23})}
W_{3})\boxtimes_{P(z_{1})}W_{4}),&\\
&(W_{1}\boxtimes_{P(z_{13})}(W_{2}\boxtimes_{P(z_{23})}
W_{3}))\boxtimes_{P(z_{3})}W_{4}),&\\
&((W_{1}\boxtimes_{P(z_{12})}W_{2})\boxtimes_{P(z_{23})}
W_{3})\boxtimes_{P(z_{3})}W_{4}),&\\
&(W_{1}\boxtimes_{P(z_{12})}W_{2})\boxtimes_{P(z_{2})}
(W_{3}\boxtimes_{P(z_{3})}W_{4})),&
\end{eqnarray*}
respectively.
\end{propo}
\pf
The proof is essentially the same as those of Corollaries \ref{prospan}
and \ref{iterspan}. 
\epfv

\begin{propo}
Let $W_1$, $W_2$, $W_3$, $W_4$ be objects
of ${\cal C}$, $w_{(1)}\in W_1$,
$w_{(2)}\in W_2$, $w_{(3)}\in W_3$, $w_{(4)}\in W_4$ and 
$z_{1}, z_{2}, z_{3}\in \C^{\times}$ such that 
$z_{12}=z_{1}-z_{2}\ne 0$, $z_{13}=z_{1}-z_{3}\ne 0$,
and $z_{23}=z_{2}-z_{3}\ne 0$.
Then:

\begin{enumerate}

\item When
$|z_{1}|>|z_{2}|>|z_{12}|+|z_{3}|>0$, 
we have 
\begin{eqnarray}\label{assoc-4-1}
\lefteqn{\overline{
\mathcal{A}_{P(z_{1}), P(z_{2})}^{P(z_{12}), P(z_{2})}}
(w_{(1)}\boxtimes_{P(z_{1})}(w_{(2)}\boxtimes_{P(z_{2})}
(w_{(3)}\boxtimes_{P(z_{3})}w_{(4)})))}\nn
&&=(w_{(1)}\boxtimes_{P(z_{12})}w_{(2)})\boxtimes_{P(z_{2})}
(w_{(3)}\boxtimes_{P(z_{3})}w_{(4)}),
\end{eqnarray}
where $\overline{
\mathcal{A}_{P(z_{1}), P(z_{2})}^{P(z_{12}), P(z_{2})}}$
is the natural extension of 
$\mathcal{A}_{P(z_{1}), P(z_{2})}^{P(z_{12}), P(z_{2})}$
to (\ref{tr-4-mod-p}).

\item When
$|z_{2}|>|z_{12}|+|z_{3}|>0$, $|z_{3}|>|z_{12}|+|z_{23}|>0$
and $|z_{23}|>|z_{12}|>0$, 
we have 
\begin{eqnarray}\label{assoc-4-2}
\lefteqn{\overline{
\mathcal{A}_{P(z_{2}), P(z_{3})}^{P(z_{23}), P(z_{3})}}
((w_{(1)}\boxtimes_{P(z_{12})}w_{(2)})\boxtimes_{P(z_{2})}
(w_{(3)}\boxtimes_{P(z_{3})}w_{(4)}))}\nn
&&
=((w_{(1)}\boxtimes_{P(z_{12})}w_{(2)})\boxtimes_{P(z_{23})}
w_{(3)})\boxtimes_{P(z_{3})}w_{(4)},
\end{eqnarray}
where $\overline{
\mathcal{A}_{P(z_{2}), P(z_{3})}^{P(z_{23}), P(z_{3})}}$
is the natural extension of 
$\mathcal{A}_{P(z_{2}), P(z_{3})}^{P(z_{23}), P(z_{3})}$
to (\ref{tr-4-mod-comp}).

\item When
$|z_{1}|>|z_{2}|>|z_{3}|>|z_{23}|>0$, $|z_{1}|>|z_{3}|+|z_{23}|>0$,
we have 
\begin{eqnarray}\label{assoc-4-3}
&{\displaystyle \overline{1_{W_{1}}\boxtimes_{P(z_{1})}
\mathcal{A}_{P(z_{2}), P(z_{3})}^{P(z_{23}), P(z_{3})}}
(w_{(1)}\boxtimes_{P(z_{1})}(w_{(2)}\boxtimes_{P(z_{2})}
(w_{(3)}\boxtimes_{P(z_{3})}w_{(4)})))}&\nn
&{\displaystyle 
=w_{(1)}\boxtimes_{P(z_{1})}((w_{(2)}\boxtimes_{P(z_{23})}
w_{(3)})\boxtimes_{P(z_{3})}w_{(4)}),}&
\end{eqnarray}
where $\overline{1_{W_{1}}\boxtimes_{P(z_{1})}
\mathcal{A}_{P(z_{2}), P(z_{3})}^{P(z_{23}), P(z_{3})}}$
is the natural extension of 
$1_{W_{1}}\boxtimes_{P(z_{1})}
\mathcal{A}_{P(z_{2}), P(z_{3})}^{P(z_{23}), P(z_{3})}$
to (\ref{tr-4-mod-p}).

\item When
$|z_{3}|>|z_{13}|>|z_{23}|>0$ and $|z_{1}|>|z_{3}|+|z_{23}|>0$,
we have 
\begin{eqnarray}\label{assoc-4-4}
\lefteqn{\overline{
\mathcal{A}_{P(z_{1}), P(z_{3})}^{P(z_{13}), P(z_{3})}}
(w_{(1)}\boxtimes_{P(z_{1})}((w_{(2)}\boxtimes_{P(z_{23})}
w_{(3)})\boxtimes_{P(z_{3})}w_{(4)}))}\nn
&&
=(w_{(1)}\boxtimes_{P(z_{13})}(w_{(2)}\boxtimes_{P(z_{23})}
w_{(3)}))\boxtimes_{P(z_{3})}w_{(4)},
\end{eqnarray}
where $\overline{
\mathcal{A}_{P(z_{1}), P(z_{3})}^{P(z_{13}), P(z_{3})}}$
is the natural extension of 
$\mathcal{A}_{P(z_{1}), P(z_{3})}^{P(z_{13}), P(z_{3})}$
to (\ref{tr-4-mod-ip}).

\item When
$|z_{3}|>|z_{13}|>|z_{23}|>|z_{12}|>0$ and  $|z_{3}|>|z_{12}|+|z_{23}|>0$, 
we have 
\begin{eqnarray}\label{assoc-4-5}
&{\displaystyle \overline{
\mathcal{A}_{P(z_{13}), P(z_{23})}^{P(z_{12}), P(z_{23})}
\boxtimes_{P(z_{3})}1_{W_{4}}}
((w_{(1)}\boxtimes_{P(z_{13})}(w_{(2)}\boxtimes_{P(z_{23})}
w_{(3)}))\boxtimes_{P(z_{3})}w_{(4)})}&\nn
&{\displaystyle 
=((w_{(1)}\boxtimes_{P(z_{12})}w_{(2)})\boxtimes_{P(z_{23})}
w_{(3)})\boxtimes_{P(z_{3})}w_{(4)},}&
\end{eqnarray}
where $\overline{
\mathcal{A}_{P(z_{13}), P(z_{23})}^{P(z_{12}), P(z_{23})}
\boxtimes_{P(z_{3})}1_{W_{4}}}$
is the natural extension of 
$\mathcal{A}_{P(z_{13}), P(z_{23})}^{P(z_{12}), P(z_{23})}
\boxtimes_{P(z_{3})}1_{W_{4}}$
to (\ref{tr-4-mod-pi}).

\end{enumerate}
\end{propo}
\pf
To prove (\ref{assoc-4-1}), 
we note that when $|z_{1}|>|z_{2}|>|z_{12}|+|z_{3}|>0$,
by Corollary \ref{t-prod-4-elts}, for 
\[
w'\in (W_{1}\boxtimes_{P(z_{12})}W_{2})\boxtimes_{P(z_{2})}
(W_{3}\boxtimes_{P(z_{3})}W_{4}),
\]
we have the absolutely convergent series
\begin{eqnarray*}
\lefteqn{\langle w', w_{(1)}\boxtimes_{P(z_{1})}(w_{(2)}\boxtimes_{P(z_{2})}
(w_{(3)}\boxtimes_{P(z_{3})}w_{(4)}))\rangle}\nn
&&=\sum_{m\in \R}
\langle w', w_{(1)}\boxtimes_{P(z_{1})}(w_{(2)}\boxtimes_{P(z_{2})}
\pi_{m}(w_{(3)}\boxtimes_{P(z_{3})}w_{(4)}))\rangle
\end{eqnarray*}
and
\begin{eqnarray*}
\lefteqn{\langle w', (w_{(1)}\boxtimes_{P(z_{12})}w_{(2)})\boxtimes_{P(z_{2})}
(w_{(3)}\boxtimes_{P(z_{3})}w_{(4)})\rangle}\nn
&&=\sum_{m\in \R}\langle w', (w_{(1)}\boxtimes_{P(z_{12})}w_{(2)})
\boxtimes_{P(z_{2})}
\pi_{m}(w_{(3)}\boxtimes_{P(z_{3})}w_{(4)})\rangle.
\end{eqnarray*}
Using the map $\left(\mathcal{A}_{P(z_{1}), P(z_{2})}^{P(z_{12}), P(z_{2})}\right)'$,
we have
\begin{eqnarray*}
\lefteqn{\left\langle w', 
\overline{\mathcal{A}_{P(z_{1}), P(z_{2})}^{P(z_{12}), P(z_{2})}}
(w_{(1)}\boxtimes_{P(z_{1})}(w_{(2)}\boxtimes_{P(z_{2})}
(w_{(3)}\boxtimes_{P(z_{3})}w_{(4)})))\right\rangle}\nn
&&=\left\langle 
\left(\mathcal{A}_{P(z_{1}), P(z_{2})}^{P(z_{12}), P(z_{2})}\right)'(w'), 
 w_{(1)}\boxtimes_{P(z_{1})}(w_{(2)}\boxtimes_{P(z_{2})}
(w_{(3)}\boxtimes_{P(z_{3})}w_{(4)}))\right\rangle\nn
&&=\sum_{m\in \R}\left\langle 
\left(\mathcal{A}_{P(z_{1}), P(z_{2})}^{P(z_{12}), P(z_{2})}\right)'(w'), 
w_{(1)}\boxtimes_{P(z_{1})}(w_{(2)}\boxtimes_{P(z_{2})}
\pi_{m}(w_{(3)}\boxtimes_{P(z_{3})}w_{(4)}))\right\rangle\nn
&&=\sum_{m\in \R}\left\langle w', 
\overline{\mathcal{A}_{P(z_{1}), P(z_{2})}^{P(z_{12}), P(z_{2})}}(
w_{(1)}\boxtimes_{P(z_{1})}(w_{(2)}\boxtimes_{P(z_{2})}
\pi_{m}(w_{(3)}\boxtimes_{P(z_{3})}w_{(4)})))\right\rangle\nn
&&=\sum_{m\in \R}\left\langle w', 
(w_{(1)}\boxtimes_{P(z_{12})}w_{(2)})\boxtimes_{P(z_{2})}
\pi_{m}(w_{(3)}\boxtimes_{P(z_{3})}w_{(4)})\right\rangle\nn
&&=\langle w', (w_{(1)}\boxtimes_{P(z_{12})}w_{(2)})\boxtimes_{P(z_{2})}
(w_{(3)}\boxtimes_{P(z_{3})}w_{(4)})\rangle.
\end{eqnarray*}
Since $w'$ is arbitrary, we obtain (\ref{assoc-4-1}).

The equalities (\ref{assoc-4-2})--(\ref{assoc-4-5})
are proved similarly.
\epfv

\subsection{The data of the braided monoidal category structure}

We choose the tensor product bifunctor of the braided tensor category
that we are constructing to be the bifunctor $\boxtimes_{P(1)}$, and
we shall write it simply as:
\[
\boxtimes = \boxtimes_{P(1)}.
\]
We take the unit object to be $V$.  For any $z\in \C^{\times}$ and any
object $(W, Y_{W})$ of ${\cal C}$, we take the {\it left $P(z)$-unit
isomorphism}
\[
l_{z; W}: V\boxtimes_{P(z)}W \to W
\]
to be the unique module map {}from $V\boxtimes_{P(z)}W$ to $W$ such
that
\[
\overline{l_{z; W}}\circ \boxtimes_{P(z)}=I_{Y_{W}, 0},
\]
where $I_{Y_{W}, 0}=I_{Y_{W}, p}$ for $p\in \mathbb{Z}$ is the unique
$P(z)$-intertwining map associated to the intertwining operator
$Y_{W}$ of type ${W\choose VW}$;  the existence and uniqueness of
$l_{z; W}$ are guaranteed by the universal property of the
$P(z)$-tensor product $\boxtimes_{P(z)}$.  It is characterized
by
\[
l_{z;
W}(\mathbf{1}\boxtimes_{P(z)} w)=w
\]
for $w\in W$. (Here we do not need the natural extension $\overline{l_{z;W}}$
because $\mathbf{1}\boxtimes_{P(z)} w\in V\boxtimes_{P(z)} W$.) 
The isomorphisms $l_{W; z}$ for $W\in \ob \mathcal{C}$
give a natural isomorphism $l_{z}$ from the functor $V\boxtimes_{P(z)}
\cdot$ to the identity functor $1_{\mathcal{C}}$ of $\mathcal{C}$; the
naturality of $l_{z}$ follows immediately from the characterization of
$l_{z; W}$ and (\ref{sigma1sigma2}).  The {\it right $P(z)$-unit
isomorphism}
\[
r_{z; W}: W\boxtimes_{P(z)}V \to W
\]
is the unique module map 
{}from $W\boxtimes_{P(z)}V$ to $W$ such that 
\[
\overline{r_{z; W}}\circ \boxtimes_{P(z)}=I_{\Omega_{0}(Y_{W}), 0},
\]
where $I_{\Omega_{0}(Y_{W}), 0}=I_{\Omega_{0}(Y_{W}), p}$ for $p\in
\mathbb{Z}$ is the unique $P(z)$-intertwining map associated to the
intertwining operator $\Omega_{0}(Y_{W})$ of type ${W\choose WV}$ 
(recall (\ref{Omega_r})). It is characterized by
\[
\overline{r_{z; W}}(w\boxtimes_{P(z)}
\mathbf{1})=e^{zL(-1)} w
\]
for $w\in W$.  The isomorphisms $r_{z; W}$ for $W\in \ob \mathcal{C}$ 
give a natural isomorphism $r_{z}$ from the functor $\cdot \boxtimes_{P(z)} V$
to the identity functor $1_{\mathcal{C}}$. In particular, we have the
left unit isomorphism
\[
l=l_{1}: V\boxtimes \cdot \to 1_{\mathcal{C}}
\]
and the right unit isomorphism
\[
r=r_{1}: \cdot\boxtimes V \to 1_{\mathcal{C}}.
\]

To give the braiding and associativity isomorphisms, we need
``parallel transport isomorphisms" between $P(z)$-tensor products with
different $z$.  Let $W_{1}$ and $W_{2}$ be objects of $\mathcal{C}$,
let $z_{1}, z_{2}\in \C^{\times}$, and let $\gamma$ be a path in
$\mathbb{C}^{\times}$ {}from $z_{1}$ to $z_{2}$.  Let $\mathcal{Y}$ be
the logarithmic intertwining operator associated to the
$P(z_{2})$-tensor product $W_{1}\boxtimes_{P(z_{2})}W_{2}$ and
$l(z_{1})$ the value of the logarithm of $z_{1}$ determined uniquely
by $\log z_{2}$ and the path $\gamma$. Then we have a
$P(z_{1})$-intertwining map $I$ defined by
\[
I(w_{(1)}\otimes w_{(2)})=\mathcal{Y}(w_{(1)}, e^{l(z_{1})})w_{(2)}
\]
for $w_{(1)}\in W_{1}$ and $w_{(2)}\in W_{2}$. The {\it parallel transport 
isomorphism
\[
\mathcal{T}_{\gamma}=\mathcal{T}_{\gamma;W_1,W_2}: W_{1}\boxtimes_{P(z_{1})}W_{2}
\to W_{1}\boxtimes_{P(z_{2})}W_{2}
\]
associated to $\gamma$} is defined to be the unique module map such
that
\[
I=\overline{\mathcal{T}_{\gamma}}\circ \boxtimes_{P(z_{1})},
\]
where $\overline{\mathcal{T}_{\gamma}}$ is the natural extension of
$\mathcal{T}_{\gamma}$ to the formal completion
$\overline{W_{1}\boxtimes_{P(z_{1})} W_{2}}$ of
$W_{1}\boxtimes_{P(z_{1})} W_{2}$.  As in the definition of the left
$P(z)$-unit isomorphism, the existence and uniqueness are guaranteed by
the universal property of the $P(z_{1})$-tensor product
$\boxtimes_{P(z_{1})}$.  The parallel transport isomorphism
$\mathcal{T}_{\gamma}$ is characterized by
\[
\overline{\mathcal{T}_{\gamma}}(w_{(1)}
\boxtimes_{P(z_{1})}w_{(2)})=\mathcal{Y}(w_{(1)}, x)w_{(2)}
|_{\log x=l(z_{1}),\;x^{n}=e^{nl(z_{1})}, \; n\in \R}
\]
for $w_{(1)}\in W_{1}$ and $w_{(2)}\in W_{2}$. These isomorphisms give
a natural isomorphism, denoted using the same notation
$\mathcal{T}_{\gamma}$, from $\boxtimes_{P(z_{1})}$ to
$\boxtimes_{P(z_{2})}$; the naturality of $\mathcal{T}_{\gamma}$
follows immediately from this characterization and
(\ref{sigma1sigma2}).  Since the intertwining map $I$ depends only on
the homotopy class of $\gamma$, {}from the definition we see that the
parallel transport isomorphism also depends only on the homotopy class
of $\gamma$.

For $z\in \C^{\times}$, let $I$ be the $P(z)$-intertwining map of type
${W_{2}\boxtimes_{P(-z)} W_{1} \choose W_{1}\;\; W_{2}}$ defined by
\[
I(w_{1}\otimes w_{2})= e^{zL(-1)}
(w_{(2)}\boxtimes_{P(-z)}w_{(1)})
\]
for $w_{1}\in W_{1}$ and $w_{2}\in W_{2}$.  We define a {\it
commutativity isomorphism between the $P(z)$- and $P(-z)$-tensor
products} to be the unique module map
\[
\mathcal{R}_{P(z)}=\mathcal{R}_{P(z); W_1,W_2}: W_{1}\boxtimes_{P(z)} W_{2}\to 
W_{2}\boxtimes_{P(-z)} W_{1}
\]
such that 
\[
I=\overline{\mathcal{R}_{P(z)}}\circ \boxtimes_{P(z)}.
\]
The existence and uniqueness of $\mathcal{R}_{P(z)}$ are guaranteed by
the universal property of the $P(z)$-tensor product.  By definition,
the commutativity isomorphism $\mathcal{R}_{P(z)}$ is characterized by
\[
\overline{\mathcal{R}_{P(z)}}(w_{(1)}\boxtimes_{P(z)} w_{(2)})=e^{zL(-1)}
(w_{(2)}\boxtimes_{P(-z)} w_{(1)})
\]
for $w_{(1)}\in W_{1}$, $w_{(2)}\in W_{2}$. These isomorphisms give a natural 
isomorphism, denoted using the same notation $\mathcal{R}_{P(z)}$,
\[
\mathcal{R}_{P(z)}:\boxtimes_{P(z)} \to \boxtimes_{P(-z)}\circ\sigma_{12},
\]
where $\sigma_{12}$ is the permutation functor on $\mathcal{C}\times\mathcal{C}$; 
the naturality of $\mathcal{R}_{P(z)}$
follows immediately from this characterization and (\ref{sigma1sigma2}).

Let $\gamma_{1}^{-}$ be a (``counterclockwise'') path {}from $-1$ to
$1$ in the closed upper half plane with $0$ deleted,
$\mathcal{T}_{\gamma_{1}^{-}}$ the corresponding parallel transport
isomorphism.  We define the {\it braiding isomorphism}
\[
\mathcal{R}=\mathcal{R}_{W_1;W_2}: W_{1}\boxtimes W_{2}\to
W_{2}\boxtimes W_{1}
\]
for our
braided tensor category to be 
\[
\mathcal{R}=\mathcal{T}_{\gamma_{1}^{-}}\circ \mathcal{R}_{P(1)}.
\]
This braiding isomorphism $\mathcal{R}$ can also be defined directly
as follows: Let $I$ be the $P(1)$-intertwining map of type
${W_{2}\boxtimes_{P(1)} W_{1} \choose W_{1}\;\; W_{2}}$ defined by
\[
I(w_{1}\otimes w_{2})= e^{L(-1)}
\overline{\mathcal{T}_{\gamma_{1}^{-}}} (w_{(2)}\boxtimes_{P(-1)} w_{(1)})
\]
for $w_{(1)}\in W_{1}$ and $w_{(2)}\in W_{2}$.  Then $\mathcal{R}$ is the
unique module map
\[
\mathcal{R}: W_{1}\boxtimes W_{2}\to 
W_{2}\boxtimes W_{1}
\]
such that 
\[
I=\overline{\mathcal{R}}\circ \boxtimes.
\]
It is characterized by
\[
\overline{\mathcal{R}}(w_{(1)}\boxtimes_{P(1)} w_{(2)})=e^{L(-1)}
\overline{\mathcal{T}_{\gamma_{1}^{-}}}
(w_{(2)}\boxtimes_{P(-1)} w_{(1)})
\]
for $w_{(1)}\in W_{1}$, $w_{(2)}\in W_{2}$. We have a natural isomorphism
$\mathcal{R}$ from $\boxtimes$ to $\boxtimes\circ\sigma_{12}$.

Let $z_{1}$ and $z_{2}$ be complex numbers satisfying
$|z_{1}|>|z_{2}|>|z_{1}-z_{2}|>0$ and let $W_{1}$, $W_{2}$ and $W_{3}$
be objects of $\mathcal{C}$.  Then the associativity isomorphism
(corresponding to the indicated geometric data)
\[
\alpha_{P(z_{1}), P(z_{2})}^{P(z_{1}-z_{2}), P(z_{2})}:
(W_{1}\boxtimes_{P(z_{1}-z_{2})}W_{2})\boxtimes_{P(z_{2})}W_{3}
\to W_{1}\boxtimes_{P(z_{1})}(W_{2}\boxtimes_{P(z_{2})}W_{3})
\]
and its inverse
\[
\A_{P(z_{1}), P(z_{2})}^{P(z_{1}-z_{2}), P(z_{2})}:
W_{1}\boxtimes_{P(z_{1})}(W_{2}\boxtimes_{P(z_{2})}W_{3})
\to (W_{1}\boxtimes_{P(z_{1}-z_{2})}W_{2})\boxtimes_{P(z_{2})}W_{3}
\]
have been constructed in Section 10 and are determined uniquely by
(\ref{assoc-elt-2}) and (\ref{assoc-elt-1}), respectively.

Let $z_{1}, z_{2}, z_{3}$ and $z_{4}$ be any nonzero complex numbers
(with no restrictions other than being nonzero; they may all be equal,
for example) and let $W_{1}$, $W_{2}$ and $W_{3}$ be objects of ${\cal
C}$.  In this generality, we define a natural associativity
isomorphism
\[
\A_{P(z_{1}), P(z_{2})}^{P(z_{4}), P(z_{3})}:
W_{1}\boxtimes_{P(z_{1})}(W_{2}\boxtimes_{P(z_{2})}W_{3})
\to (W_{1}\boxtimes_{P(z_{4})}W_{2})\boxtimes_{P(z_{3})}W_{3}
\]
using the already-constructed associativity and parallel transport
isomorphisms as follows: Let $\zeta_{1}$ and $\zeta_{2}$ be nonzero
complex numbers satisfying
$|\zeta_{1}|>|\zeta_{2}|>|\zeta_{1}-\zeta_{2}|>0$.  Let $\gamma_{1}$
and $\gamma_{2}$, be paths {}from $z_{1}$ and $z_{2}$ to $\zeta_{1}$
and $\zeta_{2}$, respectively, in the complex plane with a cut along
the positive real line, and let $\gamma_{3}$ and $\gamma_{4}$ be paths
{}from $\zeta_{2}$ and $\zeta_{1}-\zeta_{2}$ to $z_{3}$ and $z_{4}$,
respectively, also in the complex plane with a cut along the positive
real line.  Then we define
\[
\A_{P(z_{1}), P(z_{2})}^{P(z_{4}), P(z_{3})}
=\mathcal{T}_{\gamma_{3}}\circ (\mathcal{T}_{\gamma_{4}}
\boxtimes_{P(\zeta_{2})} 1_{W_{3}})\circ
\mathcal{A}^{P(\zeta_{1}-\zeta_{2}), P(\zeta_{2})}_{P(\zeta_{1}), P(\zeta_{2})}\circ
(1_{W_{1}} \boxtimes_{P(\zeta_{1})}
\mathcal{T}_{\gamma_{2}})\circ \mathcal{T}_{\gamma_{1}},
\]
that is, $\A_{P(z_{1}), P(z_{2})}^{P(z_{4}), P(z_{3})}$ 
is given by the commutative diagram
\[
\begin{CD}
W_{1}\boxtimes_{P(\zeta_{1})} (W_{2}
\boxtimes_{P(\zeta_{2})} W_{3})
@>\mathcal{A}_{P(\zeta_{1}), P(\zeta_{2})}^{P(\zeta_{1}-z_{2}), P(\zeta_{2})}>>
(W_{1}\boxtimes_{P(\zeta_{1}-\zeta_{2})} W_{2})
\boxtimes_{P(\zeta_{2})} W_{3}\\
@A(1_{W_{1}} \boxtimes_{P(\zeta_{1})}
\mathcal{T}_{\gamma_{2}})\circ \mathcal{T}_{\gamma_{1}}AA
@VV\mathcal{T}_{\gamma_{3}}\circ (\mathcal{T}_{\gamma_{4}}
\boxtimes_{P(\zeta_{2})} 1_{W_{3}})V\\
W_{1}\boxtimes_{P(z_{1})} (W_{2}
\boxtimes_{P(z_{2})} W_{3}) @>\mathcal{A}_{P(z_{1}), P(z_{2})}^{P(z_{4}), P(z_{3})}>>
(W_{1}\boxtimes_{P(z_{4})} W_{2})
\boxtimes_{P(z_{3})} W_{3}
\end{CD}
\]
The inverse of $\A_{P(z_{1}), P(z_{2})}^{P(z_{4}), P(z_{3})}$ is
denoted $\alpha_{P(z_{1}), P(z_{2})}^{P(z_{4}), P(z_{3})}$.  These
isomorphisms certainly generalize the previously-constructed ones
(when the appropriate inequalities hold).

In particular, when $z_{1}=z_{2}=z_{3}=z_{4}=1$, we call the
corresponding natural associativity isomorphism
\[
\A_{P(1), P(1)}^{P(1), P(1)}:
W_{1}\boxtimes (W_{2}\boxtimes W_{3})\to
(W_{1}\boxtimes W_{2})\boxtimes W_{3}
\]
{\it the associativity isomorphism (for the braided tensor category
structure)} and denote it simply by $\A$. Its inverse is denoted by
$\alpha$.  Because of the importance of this special case, we rewrite
the definition of $\A$ explicitly. Let $r_{1}$ and $r_{2}$ be real
numbers satisfying $r_{1}>r_{2}>r_{1}-r_{2}> 0$.  Let $\gamma_{1}$
and $\gamma_{2}$ be paths in $(0, \infty)$ {}from $1$ to $r_{1}$ and
to $r_{2}$, respectively, and let $\gamma_{3}$ and $\gamma_{4}$ be
paths in $(0, \infty)$ {}from $r_{2}$ and {}from $r_{1}-r_{2}$ to $1$,
respectively.  Then
\[
\mathcal{A}=\mathcal{T}_{\gamma_{3}}\circ (\mathcal{T}_{\gamma_{4}}
\boxtimes_{P(r_{2})} 1_{W_{3}})\circ
\mathcal{A}^{P(r_{1}-r_{2}), P(r_{2})}_{P(r_{1}), P(r_{2})}\circ
(1_{W_{1}} \boxtimes_{P(r_{1})}
\mathcal{T}_{\gamma_{2}})\circ \mathcal{T}_{\gamma_{1}},
\]
that is, $\mathcal{A}$ is given by the commutative diagram
\[
\begin{CD}
W_{1}\boxtimes_{P(r_{1})} (W_{2}
\boxtimes_{P(r_{2})} W_{3})
@>\mathcal{A}_{P(r_{1}), P(r_{2})}^{P(r_{1}-r_{2}), P(r_{2})}>>
(W_{1}\boxtimes_{P(r_{1}-r_{2})} W_{2})
\boxtimes_{P(r_{2})} W_{3}\\
@A(1_{W_{1}} \boxtimes_{P(r_{1})}
\mathcal{T}_{\gamma_{2}})\circ \mathcal{T}_{\gamma_{1}}AA
@VV\mathcal{T}_{\gamma_{3}}\circ (\mathcal{T}_{\gamma_{4}}
\boxtimes_{P(r_{2})} 1_{W_{3}})V\\
W_{1}\boxtimes (W_{2}
\boxtimes W_{3}) @>\mathcal{A}>>
(W_{1}\boxtimes W_{2})
\boxtimes W_{3}.
\end{CD}
\]

\begin{rema}\label{notensorprodelts}
{\rm It is important to note that in this case, or more generally,
whenever $z_1=z_2$ or $z_3=z_4$ (or both), the corresponding tensor
products of {\it elements} fail to exist; for instance, the symbol
$w_{(1)}\boxtimes_{P(1)} (w_{(2)} \boxtimes_{P(1)} w_{(3)})$ has no
meaning.  Because of this, it has been necessary for us to develop our
whole theory for general nonzero complex numbers (and this in turn has
required the theory of the logarithmic operator product expansion and
so on).}
\end{rema}

\subsection{Actions of the associativity and
commutativity isomorphisms on tensor products of elements}

\begin{propo}
For any $z_{1}, z_{2}\in \C^{\times}$ such that 
$z_{1}\ne z_{2}$ but $|z_{1}|=|z_{2}|=|z_{1}-z_{2}|$,
we have 
\begin{equation}\label{assoc-general-z}
\overline{\A_{P(z_{1}), P(z_{2})}^{P(z_{1}-z_{2}), P(z_{2})}}
(w_{(1)}\boxtimes_{P(z_{1})}(w_{(2)}\boxtimes_{P(z_{2})}w_{(3)}))
=(w_{(1)}\boxtimes_{P(z_{1}-z_{2})}w_{(2)})\boxtimes_{P(z_{2})}w_{(3)}
\end{equation}
for $w_{(1)}\in W_{1}$, $w_{(2)}\in W_{2}$ and $w_{(3)}\in W_{3}$,
where
\[
\overline{\A_{P(z_{1}), P(z_{2})}^{P(z_{1}-z_{2}), P(z_{2})}}:
\overline{W_{1}\boxtimes_{P(z_{1})}(W_{2}\boxtimes_{P(z_{2})}W_{3})}
\to \overline{(W_{1}\boxtimes_{P(z_{1}-z_{2})}W_{2})
\boxtimes_{P(z_{2})}W_{3}}
\]
is the natural extension of
$\A_{P(z_{1}), P(z_{2})}^{P(z_{1}-z_{2}), P(z_{2})}$.
\end{propo}
\pf We need only prove the case that $w_{(1)}$ and $w_{(2)}$ are
homogeneous with respect to the generalized-weight grading.  So we now
assume that they are homogeneous.

We can always find $\epsilon_{1}\in \C$
such that 
\begin{eqnarray}
|z_{1}+\epsilon_{1}|&>&|\epsilon_{1}|,\label{assoc-general-z-0}\\
|z_{1}+\epsilon_{1}|&>&|z_{2}|
>|(z_{1}+\epsilon_{1})-z_{2}|>0.\label{assoc-general-z-1}
\end{eqnarray}
Let 
$\Y_{1}=\Y_{\boxtimes_{P(z_{1})}, 0}$ and 
$\Y_{2}=\Y_{\boxtimes_{P(z_{2})}, 0}$
be intertwining operators of types 
\[
{W_{1}\boxtimes_{P(z_{1})}(W_{2}\boxtimes_{P(z_{2})}W_{3})\choose 
W_{1}\;\; W_{2}\boxtimes_{P(z_{2})}W_{3}}
\]
and 
\[
{W_{2}\boxtimes_{P(z_{2})}W_{3}\choose W_{2}\;\; W_{3}},
\]
respectively, corresponding to the intertwining maps 
$\boxtimes_{P(z_{1})}$ and $\boxtimes_{P(z_{2})}$, respectively. 
Then the series 
\[
\langle w', \Y_{1}(\pi_{m}(e^{-\epsilon_{1}L(-1)}w_{(1)}), z_{1}+\epsilon_{1})
\Y_{2}(\pi_{n}(w_{(2)}), z_{2})w_{(3)}\rangle
\]
is absolutely convergent for $m, n\in \R$ and 
\[
w'\in 
(W_{1}\boxtimes_{P(z_{1})}(W_{2}\boxtimes_{P(z_{2})}W_{3}))'
\]
 and the sums of these series define elements
\[
\Y_{1}(\pi_{m}(e^{-\epsilon_{1}L(-1)}w_{(1)}), z_{1}+\epsilon_{1})
\Y_{2}(w_{(2)}, z_{2})w_{(3)}
\in \overline{W_{1}\boxtimes_{P(z_{1})}(W_{2}\boxtimes_{P(z_{2})}W_{3})}.
\]

By the definition of the parallel transport isomorphism, for any path
$\gamma$ {}from $z_{1}+\epsilon_{1}$ to $z_{1}$ in the complex plane
with a cut along the nonnegative real line, we have
\begin{eqnarray}\label{assoc-general-z-2}
\lefteqn{\overline{\mathcal{T}_{\gamma}}
(\pi_{m}(e^{-\epsilon_{1}L(-1)}w_{(1)}) \boxtimes_{P(z_{1}+\epsilon_{1})}
(w_{(2)}\boxtimes_{P(z_{2})}
w_{(3)}))}\nn
&&=\Y_{1}(\pi_{m}(e^{-\epsilon_{1}L(-1)}w_{(1)}), z_{1}+\epsilon_{1})
\Y_{2}(w_{(2)}, z_{2})w_{(3)}.
\end{eqnarray}
By definition, we know that 
\[
\mathcal{T}_{\gamma}^{-1}=\mathcal{T}_{\gamma^{-1}},
\]
so that (\ref{assoc-general-z-2}) can be written as 
\begin{eqnarray}\label{assoc-general-z-3}
\lefteqn{\overline{ \mathcal{T}_{\gamma_{1}^{-1}}}
(\Y_{1}(\pi_{m}(e^{-\epsilon_{1}L(-1)}w_{(1)}), z_{1}+\epsilon_{1})
\Y_{2}(w_{(2)}, z_{2})w_{(3)})}\nn
&&=\pi_{m}(e^{-\epsilon_{1}L(-1)}w_{(1)})\boxtimes_{P(z_{1}+\epsilon_{1})}
(w_{(2)}\boxtimes_{P(z_{2})}
w_{(3)}).
\end{eqnarray}

Since (\ref{assoc-general-z-1}) holds, by (\ref{assoc-elt-1}), we have
\begin{eqnarray}\label{assoc-general-z-4}
\lefteqn{\overline{\mathcal{A}_{P(z_{1}+\epsilon_{1}), P(z_{2})}
^{P((z_{1}+\epsilon_{1})-z_{1}), P(z_{2})}}
(\pi_{m}(e^{-\epsilon_{1}L(-1)}w_{(1)})\boxtimes_{P(z_{1}+\epsilon_{1})}
(w_{(2)}\boxtimes_{P(z_{2})}
w_{(3)}))}\nn
&&=(\pi_{m}(e^{-\epsilon_{1}L(-1)}w_{(1)})
\boxtimes_{P((z_{1}+\epsilon_{1})-z_{2})}
w_{(2)})\boxtimes_{P(z_{2})}
w_{(3)}.\quad\quad\quad
\end{eqnarray}

Let $\Y_{3}=\Y_{\boxtimes_{P(z_{2})}, 0}$ and
$\Y_{4}=\Y_{\boxtimes_{P(z_{1}-z_{2})}, 0}$ be intertwining operators
of types
\[
{(W_{1}\boxtimes_{P(z_{1}-z_{2})}W_{2})\boxtimes_{P(z_{2})}W_{3}\choose 
(W_{1}\boxtimes_{P(z_{1}-z_{2})}W_{2})\;\;W_{3}}
\]
and 
\[
{W_{1}\boxtimes_{P(z_{1}-z_{2})}W_{2}\choose W_{1}\;\; W_{2}},
\]
respectively, corresponding to the intertwining maps 
$\boxtimes_{P(z_{2})}$ and $\boxtimes_{P(z_{1}-z_{2})}$, respectively. 
Then the series 
\[
\langle \tilde{w}', \Y_{3}(\Y_{4}(\pi_{m}(e^{-\epsilon_{1}L(-1)}w_{(1)}), 
(z_{1}+\epsilon_{1})-z_{2})
w_{(2)}, z_{2})w_{(3)}\rangle
\]
is absolutely convergent for $m, n\in \R$ and 
\[
\tilde{w}'\in 
((W_{1}\boxtimes_{P(z_{1}-z_{2})}W_{2})\boxtimes_{P(z_{2})}W_{3})',
\]
 and the sums of these series define elements
\begin{eqnarray*}
&\Y_{3}(\Y_{4}(\pi_{m}(e^{-\epsilon_{1}L(-1)}w_{(1)}), 
(z_{1}+\epsilon_{1})-z_{2})
w_{(2)}), z_{2})w_{(3)}&\nn
& \in \overline{(W_{1}\boxtimes_{P(z_{1}-z_{2})}W_{2})
\boxtimes_{P(z_{2})}W_{3}}.&
\end{eqnarray*}

We can always choose $\gamma$ such that the path $\gamma-z_{2}$ {}from
$(z_{1}+\epsilon_{1}) -z_{2}$ to $z_{1}-z_{2}$ is also in the complex
plane with a cut along the nonnegative real line. Then by the
definition of the parallel transport isomorphism, we have
\begin{eqnarray}\label{assoc-general-z-5}
\lefteqn{\overline{\mathcal{T}_{\gamma-z_{2}}}
((\pi_{m}(e^{-\epsilon_{1}L(-1)}w_{(1)})\boxtimes_{P((z_{1}+\epsilon_{1})
-z_{2})}w_{(2)})\boxtimes_{P(z_{2})}
w_{(3)})}\nn
&&=\Y_{3}(\Y_{4}(\pi_{m}(e^{-\epsilon_{1}L(-1)}w_{(1)}), (z_{1}+\epsilon_{1})
-z_{2})w_{(2)}, z_{2})w_{(3)}.\nn
\end{eqnarray}

Combining (\ref{assoc-general-z-3})--(\ref{assoc-general-z-5}) and 
using the definition of the associativity isomorphism
$\A_{P(z_{1}), P(z_{2})}^{P(z_{1}-z_{2}), P(z_{2})}$,
we obtain
\begin{eqnarray}\label{assoc-general-z-6}
\lefteqn{\overline{\A_{P(z_{1}), P(z_{2})}
^{P(z_{1}-z_{2}), P(z_{2})}}
(\Y_{1}(\pi_{m}(e^{-\epsilon_{1}L(-1)}w_{(1)}), z_{1}+\epsilon_{1})
\Y_{2}(w_{(2)}, z_{2})w_{(3)})}\nn
&&=\Y_{3}(\Y_{4}(\pi_{m}(e^{-\epsilon_{1}L(-1)}w_{(1)}), (z_{1}+\epsilon_{1})
-z_{2})w_{(2)}), z_{2})w_{(3)}.\nn
\end{eqnarray}

Since (\ref{assoc-general-z-1}) holds, both the series
\[
\sum_{m\in \R}\langle w',
\Y_{1}(\pi_{m}(e^{-\epsilon_{1}L(-1)}w_{(1)}), z_{1}+\epsilon_{1})
\Y_{2}(w_{(2)}, z_{2})w_{(3)}
\rangle
\]
and 
\[
\sum_{m\in \R}\langle \tilde{w}', 
\Y_{3}(\Y_{4}(\pi_{m}(e^{-\epsilon_{1}L(-1)}w_{(1)}), (z_{1}+\epsilon_{1})
-z_{2})
w_{(2)}, z_{2})w_{(3)}\rangle
\]
are absolutely convergent for 
\[
w'\in 
(W_{1}\boxtimes_{P(z_{1})}(W_{2}\boxtimes_{P(z_{2})}W_{3}))'
\]
and 
\[
\tilde{w}'\in 
((W_{1}\boxtimes_{P(z_{1}-z_{2})}W_{2})\boxtimes_{P(z_{2})}W_{3})'.
\]

We know that
\[
\langle \tilde{w}', 
\Y_{1}(w_{(1)}, z_{1}+\epsilon_{1})
\Y_{2}(w_{(2)}, z_{2})w_{(3)}
\rangle
\]
and 
\[
\langle \tilde{w}', 
\Y_{3}(\Y_{4}(w_{(1)}, (z_{1}+\epsilon_{1})
-z_{2})
w_{(2)}, z_{2})w_{(3)}\rangle
\]
are the values of single-valued analytic 
functions
\[
F(w', w_{(1)}, w_{(2)}, w_{(3)}; \zeta_{1}, \zeta_{2})
\]
and 
\[
G(\tilde{w}', w_{(1)}, w_{(2)}, w_{(3)}; \zeta_{1}, \zeta_{2})
\]
of $\zeta_{1}$ and $\zeta_{2}$ in a neighborhood of the point
$(\zeta_{1}, \zeta_{2})=(z_{1}+\epsilon_{1}, z_{2})$ which contains
the point $(\zeta_{1}, \zeta_{2})=(z_{1}, z_{2})$.  Then by the
definition of the tensor product elements
$w_{(1)}\boxtimes_{P(z_{1})}(w_{(2)}\boxtimes_{P(z_{2})}w_{(3)})$ and
$(w_{(1)}\boxtimes_{P(z_{1}-z_{2})}w_{(2)})\boxtimes_{P(z_{2})}w_{(3)}$,
we have
\begin{equation}\label{assoc-general-z-6.1}
\langle w', w_{(1)}\boxtimes_{P(z_{1})}(w_{(2)}\boxtimes_{P(z_{2})}w_{(3)})
\rangle
=F(w', w_{(1)}, w_{(2)}, w_{(3)}; z_{1}, z_{2})
\end{equation}
and 
\begin{equation}\label{assoc-general-z-6.2}
\langle \tilde{w}', 
(w_{(1)}\boxtimes_{P(z_{1}-z_{2})}w_{(2)})\boxtimes_{P(z_{2})}w_{(3)}
\rangle
=G(\tilde{w}', w_{(1)}, w_{(2)}, w_{(3)}; z_{1}, z_{2}).
\end{equation}
On the other hand, since $F(w', w_{(1)}, w_{(2)}, w_{(3)}; \zeta_{1},
\zeta_{2})$ and $G(\tilde{w}', w_{(1)}, w_{(2)}, w_{(3)}; \zeta_{1},
\zeta_{2})$ are analytic extensions of matrix elements of products and
iterates of intertwining maps, properties of these products and
iterates also hold for these functions if they still make sense.  In
particular, they satisfy the $L(-1)$-derivative property:
\begin{eqnarray}
\frac{\partial}{\partial \zeta_{1}}
F(w', w_{(1)}, w_{(2)}, w_{(3)}; \zeta_{1}, \zeta_{2})
&=&F(w', L(-1)w_{(1)}, w_{(2)}, w_{(3)}; \zeta_{1}, \zeta_{2}),
\label{assoc-general-z-7}\\
\frac{\partial}{\partial \zeta_{1}}
G(\tilde{w}', w_{(1)}, w_{(2)}, w_{(3)}; \zeta_{1}, \zeta_{2})
&=&G(\tilde{w}', L(-1)w_{(1)}, w_{(2)}, w_{(3)}; \zeta_{1}, \zeta_{2}),
\label{assoc-general-z-9}
\end{eqnarray}
{}From the Taylor theorem (which applies since 
(\ref{assoc-general-z-0}) holds) and 
(\ref{assoc-general-z-7})--(\ref{assoc-general-z-9}), we have 
\begin{eqnarray}
F(w', w_{(1)}, w_{(2)}, w_{(3)}; z_{1}, z_{2})
&=&\sum_{m\in \R}F(w', \pi_{m}(e^{-\epsilon_{1}L(-1)})w_{(1)}, 
w_{(2)}, w_{(3)}; 
z_{1}+\epsilon_{1}, z_{2}),\quad\quad\label{assoc-general-z-11}\\
G(\tilde{w}', w_{(1)}, w_{(2)}, w_{(3)}; z_{1}, z_{2})
&=&\sum_{m\in \R}G(\tilde{w}', \pi_{m}(e^{-\epsilon_{1}L(-1)})w_{(1)}, 
w_{(2)}, w_{(3)}; 
z_{1}+\epsilon_{1}, z_{2}).\quad\quad
\label{assoc-general-z-12}
\end{eqnarray}
Thus by the definitions of 
\begin{eqnarray*}
F(w', \pi_{m}(e^{-\epsilon_{1}L(-1)})w_{(1)}, 
w_{(2)}, w_{(3)}; 
z_{1}+\epsilon_{1}, z_{2}),\\
G(\tilde{w}', \pi_{m}(e^{-\epsilon_{1}L(-1)})w_{(1)}, 
w_{(2)}, w_{(3)}; 
z_{1}+\epsilon_{1}, z_{2}),
\end{eqnarray*}
and by (\ref{assoc-general-z-11}),
(\ref{assoc-general-z-12}), (\ref{assoc-general-z-6.1})
and (\ref{assoc-general-z-6.2}), we obtain
\begin{eqnarray}
\lefteqn{\sum_{m\in \R}\langle w',
\Y_{1}(\pi_{m}(e^{-\epsilon_{1}L(-1)}w_{(1)}), z_{1}+\epsilon_{1})
\Y_{2}(w_{(2)}, z_{2})w_{(3)}
\rangle}\nn
&&\quad\quad
=\langle w', w_{(1)}\boxtimes_{P(z_{1})}(w_{(2)}\boxtimes_{P(z_{2})}w_{(3)})
\rangle
\label{assoc-general-z-13}
\end{eqnarray}
and
\begin{eqnarray}
\lefteqn{\sum_{m\in \R}\langle \tilde{w}', 
\Y_{3}(\Y_{4}(\pi_{m}(e^{-\epsilon_{1}L(-1)}w_{(1)}), (z_{1}+\epsilon_{1})
-z_{2})w_{(2)}, z_{2})w_{(3)}\rangle}\nn
&&\quad\quad\quad
=\langle \tilde{w}', 
(w_{(1)}\boxtimes_{P(z_{1}-z_{2})}w_{(2)})\boxtimes_{P(z_{2})}w_{(3)}
\rangle.
\label{assoc-general-z-14}
\end{eqnarray}
Since $w'$ and $\tilde{w}'$ are arbitrary, (\ref{assoc-general-z-13}) and
(\ref{assoc-general-z-14}) gives
\begin{eqnarray}
\lefteqn{\sum_{m\in \R}
\Y_{1}(\pi_{m}(e^{-\epsilon_{1}L(-1)}w_{(1)}), z_{1}+\epsilon_{1})
\Y_{2}(w_{(2)}, z_{2})w_{(3)}}\nn
&&\quad\quad
=w_{(1)}\boxtimes_{P(z_{1})}(w_{(2)}\boxtimes_{P(z_{2})}w_{(3)})
\quad\quad\quad\quad\quad\quad\quad\quad\quad
\label{assoc-general-z-15}
\end{eqnarray}
and
\begin{eqnarray}
\lefteqn{\sum_{m\in \R} 
\Y_{3}(\Y_{4}(\pi_{m}(e^{-\epsilon_{1}L(-1)}w_{(1)}), (z_{1}+\epsilon_{1})
-z_{2})w_{(2)}, z_{2})w_{(3)}}\nn
&&\quad\quad\quad
=(w_{(1)}\boxtimes_{P(z_{1}-z_{2})}w_{(2)})\boxtimes_{P(z_{2})}w_{(3)}.
\quad\quad\quad\quad\quad\quad\quad\quad\quad
\label{assoc-general-z-16}
\end{eqnarray}

Taking the sum $\sum_{m\in \R}$ on both sides of
(\ref{assoc-general-z-6}) and then using (\ref{assoc-general-z-15})
and (\ref{assoc-general-z-16}), we obtain (\ref{assoc-general-z}).
\epfv

We also have:

\begin{propo}
Let $z_{1}, z_{2}$ be nonzero complex numbers such that $z_{1}\ne
z_{2}$ but $|z_{1}|=|z_{2}|=|z_{1}-z_{2}|$. Let $\gamma$ be a path
{}from $z_{2}$ to $z_{1}$ in the complex plane with a cut along the
nonnegative real line. Then we have
\begin{equation}\label{commu-1}
\overline{\mathcal{T}_{\gamma}\circ 
(\mathcal{R}_{P(z_{1}-z_{2})}\boxtimes_{P(z_{2})}1_{W_{3}})}
((w_{(1)}\boxtimes_{P(z_{1}-z_{2})}w_{(2)})\boxtimes_{P(z_{2})}w_{(3)})
=(w_{(2)}\boxtimes_{P(z_{2}-z_{1})}w_{(1)})\boxtimes_{P(z_{1})}w_{(3)}
\end{equation}
for $w_{(1)}\in W_{1}$, $w_{(2)}\in W_{2}$ and $w_{(3)}\in W_{3}$.
\end{propo}
\pf
We can find $\epsilon$ such that $|z_{2}|>|\epsilon|$,
\[
|z_{2}+\epsilon|>|z_{1}-z_{2}|>0.
\]
Let 
$\Y_{1}=\Y_{\boxtimes_{P(z_{2})}, 0}$, 
$\tilde{\Y}_{1}=\Y_{\boxtimes_{P(z_{2})}, 0}$ and
$\Y_{2}=\Y_{\boxtimes_{P(z_{1}-z_{2})}, 0}$
be intertwining operators of types 
\[
{(W_{1}\boxtimes_{P(z_{1}-z_{2})}W_{2})\boxtimes_{P(z_{2})}W_{3}\choose 
(W_{1}\boxtimes_{P(z_{1}-z_{2})}W_{2})\;\;W_{3}},
\]
\[
{(W_{2}\boxtimes_{P(z_{2}-z_{1})}W_{1})\boxtimes_{P(z_{2})}W_{3}\choose 
(W_{2}\boxtimes_{P(z_{2}-z_{1})}W_{1})\;\;W_{3}},
\]
and 
\[
{W_{1}\boxtimes_{P(z_{1}-z_{2})}W_{2}\choose W_{1}\;\; W_{2}},
\]
respectively, corresponding to the intertwining maps
$\boxtimes_{P(z_{2})}$, $\boxtimes_{P(z_{2})}$ and
$\boxtimes_{P(z_{1})}$, respectively.  Note that since
$|z_{2}+\epsilon|>|z_{1}-z_{2}|>0$,
\begin{eqnarray}\label{commu-1-0-0}
\lefteqn{\Y_{1}(w_{(1)}\boxtimes_{P(z_{1}-z_{2})}w_{(2)}, 
z_{2}+\epsilon)w_{(3)}}\nn
&&=\sum_{m, n\in \R}\sum_{k=1}^{N}
(\pi_{n}(w_{(1)}\boxtimes_{P(z_{1}-z_{2})}w_{(2)}))^{\Y_{1}}_{m;k}
w_{(3)}e^{(-m-1)\log (z_{2}+\epsilon)}(\log (z_{2}+\epsilon))^{k}.
\end{eqnarray}
For $n\in \R$, 
\begin{eqnarray}\label{commu-1-0-1}
&{\displaystyle \sum_{m}\sum_{k=1}^{N}
(\pi_{n}(w_{(1)}\boxtimes_{P(z_{1}-z_{2})}w_{(2)}))^{\Y_{1}}_{m;k}
w_{(3)}e^{(-m-1)\log z_{2}}(\log z_{2})^{k}}&\nn
&{\displaystyle =(\pi_{n}(w_{(1)}\boxtimes_{P(z_{1}-z_{2})}w_{(2)}))
\boxtimes_{P(z_{2})}w_{(3)}}.&
\end{eqnarray}
But by the definition of 
$\mathcal{R}_{P(z_{1}-z_{2})}\boxtimes_{P(z_{2})}
1_{W_{3}}$,
\begin{eqnarray}\label{commu-1-0-2}
\lefteqn{\overline{(\mathcal{R}_{P(z_{1}-z_{2})}
\boxtimes_{P(z_{2})}1_{W_{3}})}
((\pi_{n}(w_{(1)}\boxtimes_{P(z_{1}-z_{2})}w_{(2)}))
\boxtimes_{P(z_{2})}w_{(3)})}\nn
&&\quad =\mathcal{R}_{P(z_{1}-z_{2})}
(\pi_{n}(w_{(1)}\boxtimes_{P(z_{1}-z_{2})}w_{(2)}))
\boxtimes_{P(z_{2})}w_{(3)}\nn
&&\quad =\pi_{n}(\mathcal{R}_{P(z_{1}-z_{2})}
(w_{(1)}\boxtimes_{P(z_{1}-z_{2})}w_{(2)}))
\boxtimes_{P(z_{2})}w_{(3)}\quad\quad
\end{eqnarray}
for $n\in \R$. {}From (\ref{commu-1-0-1}), (\ref{commu-1-0-2}) and 
the definitions of $\Y_{1}$ and $\tilde{\Y}_{1}$,
we obtain
\begin{eqnarray}\label{commu-1-0-3}
\lefteqn{\overline{(\mathcal{R}_{P(z_{1}-z_{2})}
\boxtimes_{P(z_{2})}1_{W_{3}})}
((\pi_{n}(w_{(1)}\boxtimes_{P(z_{1}-z_{2})}w_{(2)}))^{\Y_{1}}_{m;k}
w_{(3)})}\nn
&&=(\pi_{n}(\mathcal{R}_{P(z_{1}-z_{2})}
(w_{(1)}\boxtimes_{P(z_{1}-z_{2})}w_{(2)})))^{\tilde{\Y}_{1}}_{m;k}
w_{(3)}\nn
&&=(\pi_{n}(e^{(z_{1}-z_{2})L(-1)}
(w_{(2)}\boxtimes_{P(z_{2}-z_{1})}w_{(1)})))^{\tilde{\Y}_{1}}_{m;k}
w_{(3)}
\end{eqnarray}
for $m\in \R$ and $k=1, \dots, N$.
Using (\ref{commu-1-0-0}) and (\ref{commu-1-0-3}),
we obtain
\begin{eqnarray}\label{commu-1-1}
\lefteqn{\overline{(\mathcal{R}_{P(z_{1}-z_{2})}\boxtimes_{P(z_{2})}1_{W_{3}})}
(\Y_{1}(w_{(1)}\boxtimes_{P(z_{1}-z_{2})}w_{(2)}, z_{2}+\epsilon)w_{(3)})}\nn
&&=\overline{(\mathcal{R}_{P(z_{1}-z_{2})}\boxtimes_{P(z_{2})}1_{W_{3}})}
\nn
&&\quad\quad\left(\sum_{m, n\in \R}\sum_{k=1}^{N}
(\pi_{n}(w_{(1)}\boxtimes_{P(z_{1}-z_{2})}w_{(2)}))^{\Y_{1}}_{m;k}
w_{(3)}
e^{(-m-1)\log (z_{2}+\epsilon)}(\log (z_{2}+\epsilon))^{k}\right)\nn
&&=\sum_{m, n\in \R}\sum_{k=1}^{N}
\overline{(\mathcal{R}_{P(z_{1}-z_{2})}\boxtimes_{P(z_{2})}1_{W_{3}})}
((\pi_{n}(w_{(1)}\boxtimes_{P(z_{1}-z_{2})}w_{(2)}))^{\Y_{1}}_{m;k}
w_{(3)})\cdot\nn
&&\quad\quad\quad\quad\quad\quad\cdot
e^{(-m-1)\log (z_{2}+\epsilon)}(\log (z_{2}+\epsilon))^{k}\nn
&&=\sum_{m, n\in \R}\sum_{k=1}^{N}
(\pi_{n}(e^{(z_{1}-z_{2})L(-1)}
(w_{(2)}\boxtimes_{P(z_{2}-z_{1})}w_{(1)})))^{\tilde{\Y}_{1}}_{m;k}
w_{(3)}\cdot\nn
&&\quad\quad\quad\quad\quad\quad\cdot
e^{(-m-1)\log (z_{2}+\epsilon)}(\log (z_{2}+\epsilon))^{k}\nn
&&=\tilde{\Y}_{1}(e^{(z_{1}-z_{2})L(-1)}(w_{(2)}\boxtimes_{P(z_{2}-z_{1})}
w_{(1)}), z_{2}+\epsilon)w_{(3)}.
\end{eqnarray}

Let $\Y_{3}=\Y_{\boxtimes_{P(z_{1})}, 0}$ be the intertwining operator
of type
\[
{(W_{2}\boxtimes_{P(z_{2}-z_{1})}W_{1})\boxtimes_{P(z_{1})}W_{3}\choose 
(W_{2}\boxtimes_{P(z_{2}-z_{1})}W_{1})\;\;W_{3}}.
\]
Then by the definition of the parallel transport isomorphism,
for $m\in \R$
we have 
\begin{eqnarray}\label{commu-1-2}
\lefteqn{\overline{\mathcal{T}_{\gamma}}(
\tilde{\Y}_{1}(\pi_{m}(e^{(z_{1}-z_{2})L(-1)}(w_{(2)}\boxtimes_{P(z_{2}-z_{1})}
w_{(1)})), z_{2}+\epsilon)w_{(3)})}\nn
&&=\overline{\mathcal{T}_{\gamma}}(
\tilde{\Y}_{1}(e^{\epsilon L(-1)}
\pi_{m}(e^{(z_{1}-z_{2})L(-1)}(w_{(2)}\boxtimes_{P(z_{2}-z_{1})}
w_{(1)})), z_{2})w_{(3)})\nn
&&=\overline{\mathcal{T}_{\gamma}}(
(e^{\epsilon L(-1)}
\pi_{m}(e^{(z_{1}-z_{2})L(-1)}(w_{(2)}\boxtimes_{P(z_{2}-z_{1})}
w_{(1)})))\boxtimes_{P(z_{2})}w_{(3)})\nn
&&=\Y_{3}(e^{\epsilon L(-1)}
\pi_{m}(e^{(z_{1}-z_{2})L(-1)}(w_{(2)}\boxtimes_{P(z_{2}-z_{1})}
w_{(1)})), z_{2})w_{(3)}\nn
&&=\Y_{3}(
\pi_{m}(e^{(z_{1}-z_{2})L(-1)}(w_{(2)}\boxtimes_{P(z_{2}-z_{1})}
w_{(1)})), z_{2}+\epsilon)w_{(3)}.
\end{eqnarray}
Since $|z_{2}+\epsilon|>|z_{1}-z_{2}|>0$, the sums of both sides of
(\ref{commu-1-2}) for $m\in \R$ are absolutely convergent (in the
sense that the series obtained by paring it with elements of
$((W_{2}\boxtimes_{P(z_{2}-z_{1})}W_{1})\boxtimes_{P(z_{1})}W_{3})'$
are absolutely convergent) and we have
\begin{eqnarray}\label{commu-1-3}
\lefteqn{\sum_{m\in \R}\overline{\mathcal{T}_{\gamma}}(
\tilde{\Y}_{1}(\pi_{m}(e^{(z_{1}-z_{2})L(-1)}(w_{(2)}\boxtimes_{P(z_{2}-z_{1})}
w_{(1)})), z_{2}+\epsilon)w_{(3)})}\nn
&&=\overline{\mathcal{T}_{\gamma}}(
\tilde{\Y}_{1}(e^{(z_{1}-z_{2})L(-1)}(w_{(2)}\boxtimes_{P(z_{2}-z_{1})}
w_{(1)}), z_{2}+\epsilon)w_{(3)})
\end{eqnarray}
and
\begin{eqnarray}\label{commu-1-4}
\lefteqn{\sum_{m\in \R}\Y_{3}(
\pi_{m}(e^{(z_{1}-z_{2})L(-1)}(w_{(2)}\boxtimes_{P(z_{2}-z_{1})}
w_{(1)})), z_{2}+\epsilon)w_{(3)}}\nn
&&=\Y_{3}((w_{(2)}\boxtimes_{P(z_{2}-z_{1})}
w_{(1)}), z_{2}+\epsilon+(z_{1}-z_{2}))w_{(3)}\nn
&&=\Y_{3}((w_{(2)}\boxtimes_{P(z_{2}-z_{1})}
w_{(1)}), z_{1}+\epsilon)w_{(3)}.
\end{eqnarray}
{}From (\ref{commu-1-2})--(\ref{commu-1-4}), we obtain
\begin{eqnarray}\label{commu-1-5}
\lefteqn{\overline{\mathcal{T}_{\gamma}}(
\tilde{\Y}_{1}(e^{(z_{1}-z_{2})L(-1)}(w_{(2)}\boxtimes_{P(z_{2}-z_{1})}
w_{(1)}), z_{2}+\epsilon)w_{(3)})}\nn
&&=\Y_{3}((w_{(2)}\boxtimes_{P(z_{2}-z_{1})}
w_{(1)}), z_{1}+\epsilon)w_{(3)}.
\end{eqnarray}
{}From (\ref{commu-1-1}) and (\ref{commu-1-5}), we obtain
\begin{eqnarray*}
\lefteqn{\overline{\mathcal{T}_{\gamma}\circ
(\mathcal{R}_{P(z_{1}-z_{2})}\boxtimes_{P(z_{2})}1_{W_{3}})}
(\Y_{1}(w_{(1)}\boxtimes_{P(z_{1}-z_{2})}w_{(2)}, z_{2}+\epsilon)w_{(3)})}\nn
&&=\overline{\mathcal{T}_{\gamma}}(
\tilde{\Y}_{1}(e^{(z_{1}-z_{2})L(-1)}(w_{(2)}\boxtimes_{P(z_{2}-z_{1})}
w_{(1)}), z_{2}+\epsilon)w_{(3)})\nn
&&=\Y_{3}((w_{(2)}\boxtimes_{P(z_{2}-z_{1})}
w_{(1)}), z_{1}+\epsilon)w_{(3)}.
\end{eqnarray*}
Then for any 
\[
w'\in ((W_{1}\boxtimes_{P(z_{2}-z_{1})}W_{2})\boxtimes_{P(z_{1})}W_{3})',
\]
we have
\begin{eqnarray*}
\lefteqn{\langle w', \overline{\mathcal{T}_{\gamma}\circ
(\mathcal{R}_{P(z_{1}-z_{2})}\boxtimes_{P(z_{2})}1_{W_{3}})}
(\Y_{1}(w_{(1)}\boxtimes_{P(z_{1}-z_{2})}w_{(2)}, z_{2}+\epsilon)w_{(3)})
\rangle}\nn
&&\quad\quad=\langle w', \Y_{3}((w_{(2)}\boxtimes_{P(z_{2}-z_{1})}
w_{(1)}), z_{1}+\epsilon)w_{(3)}\rangle,\quad\quad\quad\quad\quad\quad
\end{eqnarray*}
or equivalently,
\begin{eqnarray}\label{commu-1-7}
\lefteqn{\langle ((\mathcal{R}_{P(z_{1}-z_{2})}\boxtimes_{P(z_{2})}1_{W_{3}})'
\circ \mathcal{T}_{\gamma}')(w'), 
\Y_{1}(w_{(1)}\boxtimes_{P(z_{1}-z_{2})}w_{(2)}, z_{2}+\epsilon)w_{(3)}
\rangle}\nn
&&\quad\quad=\langle w', \Y_{3}((w_{(2)}\boxtimes_{P(z_{2}-z_{1})}
w_{(1)}), z_{1}+\epsilon)w_{(3)}\rangle,\quad\quad\quad\quad\quad\quad
\end{eqnarray}
The left- and right-hand sides of (\ref{commu-1-7}) are values at 
$(\zeta_{1}, \zeta_{2}) =(z_{1}+\epsilon, z_{2}+\epsilon)$ of some
single-valued analytic functions of $\zeta_{1}$ and $\zeta_{2}$
defined in the region
\[
\{(\zeta_{1}, \zeta_{2})\in \C^{2}\;|\;\zeta_{1}\ne 0,\;
\zeta_{2}\ne 0,\;\zeta_{1}\ne \zeta_{2},\;
0\le \arg \zeta_{1}, \arg \zeta_{2}, 
\arg (\zeta_{1}-\zeta_{2})<2\pi\}.
\]
Also, by the definition of tensor product of three elements
above, the values of these analytic functions at $(\zeta_{1}, \zeta_{2})
=(z_{1}, z_{2})$  are equal to 
\[
\langle ((\mathcal{R}_{P(z_{1}-z_{2})}\boxtimes_{P(z_{2})}1_{W_{3}})'
\circ \mathcal{T}_{\gamma}')(w'), 
(w_{(1)}\boxtimes_{P(z_{1}-z_{2})}w_{(2)})\boxtimes_{P(z_{2})}w_{(3)}
\rangle
\]
and 
\[
\langle w', (w_{(2)}\boxtimes_{P(z_{2}-z_{1})}
w_{(1)})\boxtimes_{P(z_{1})}w_{(3)}\rangle,
\]
respectively.  Thus we can take the limit $\epsilon\to 0$ on both
sides of (\ref{commu-1-7}) and obtain
\begin{eqnarray}\label{commu-1-8}
\lefteqn{\langle ((\mathcal{R}_{P(z_{1}-z_{2})}\boxtimes_{P(z_{2})}1_{W_{3}})'
\circ \mathcal{T}_{\gamma}')(w'), 
(w_{(1)}\boxtimes_{P(z_{1}-z_{2})}w_{(2)})\boxtimes_{P(z_{2})}w_{(3)}
\rangle}\nn
&&\quad\quad=\langle w', (w_{(2)}\boxtimes_{P(z_{2}-z_{1})}
w_{(1)})\boxtimes_{P(z_{1})}w_{(3)}\rangle.\quad\quad\quad\quad\quad\quad
\end{eqnarray}
Since $w'$ is arbitrary, (\ref{commu-1-8}) is equivalent to 
(\ref{commu-1}).
\epfv

We also prove:

\begin{propo}
Let $z_{1}, z_{2}$ be nonzero complex numbers such that $z_{1}\ne
z_{2}$ but $|z_{1}|=|z_{2}|=|z_{1}-z_{2}|$. Let $\gamma$ be a path
{}from $z_{2}$ to $z_{2}-z_{1}$ in the complex plane with a cut along
the nonnegative real line. Then we have
\begin{eqnarray}\label{commu-2}
\lefteqn{\overline{\mathcal{T}_{\gamma}\circ 
(1_{W_{2}}\boxtimes_{P(z_{2})}\mathcal{R}_{P(z_{1})})}
(w_{(2)}\boxtimes_{P(z_{2})}(w_{(1)}\boxtimes_{P(z_{1})}w_{(3)}))}\nn
&&=e^{z_{1}L(-1)}(w_{(2)}\boxtimes_{P(z_{2}-z_{1})}(w_{(3)}
\boxtimes_{P(-z_{1})}w_{(1)}))
\end{eqnarray}
for $w_{(1)}\in W_{1}$, $w_{(2)}\in W_{2}$ and $w_{(3)}\in W_{3}$.
\end{propo}
\pf
We can find $\epsilon$ such that $|z_{2}|>|\epsilon|$,
\[
|z_{2}+\epsilon|, |z_{2}-z_{1}+\epsilon|>|z_{1}|>0.
\]
Let 
$\Y_{1}=\Y_{\boxtimes_{P(z_{2})}, 0}$, $\tilde{\Y}_{1}
=\Y_{\boxtimes_{P(z_{2})}, 0}$ and 
$\Y_{2}=\Y_{\boxtimes_{P(z_{1})}, 0}$
be intertwining operators of types 
\[
{W_{2}\boxtimes_{P(z_{2})}(W_{1}\boxtimes_{P(z_{1})}W_{3})\choose 
W_{2}\;\;(W_{1}\boxtimes_{P(z_{1})}W_{3})},
\]
\[
{W_{2}\boxtimes_{P(z_{2})}(W_{3}\boxtimes_{P(-z_{1})}W_{1})\choose 
W_{2}\;\;(W_{3}\boxtimes_{P(-z_{1})}W_{1})}
\]
and 
\[
{W_{1}\boxtimes_{P(z_{1})}W_{3}\choose W_{1}\;\; W_{3}},
\]
respectively, corresponding to the intertwining maps
$\boxtimes_{P(z_{2})}$, $\boxtimes_{P(z_{2})}$ and
$\boxtimes_{P(z_{1})}$, respectively.  Since
$|z_{2}+\epsilon|>|z_{1}|>0$,
\begin{eqnarray}\label{commu-2-0-0}
\lefteqn{\Y_{1}(w_{(2)}, 
z_{2}+\epsilon)(w_{(1)}\boxtimes_{P(z_{1})}w_{(3)})}\nn
&&=\sum_{m, n\in \R}\sum_{k=1}^{N}
(w_{(2)})^{\Y_{1}}_{m;k}\pi_{n}(w_{(1)}\boxtimes_{P(z_{1})}
w_{(3)})e^{(-m-1)\log (z_{2}+\epsilon)}(\log (z_{2}+\epsilon))^{k}.
\end{eqnarray}
For $n\in \R$, 
\begin{eqnarray}\label{commu-2-0-1}
&{\displaystyle \sum_{m\in\R}\sum_{k=1}^{N}
(w_{(2)})^{\Y_{1}}_{m;k}\pi_{n}(w_{(1)}\boxtimes_{P(z_{1})}
w_{(3)})e^{(-m-1)\log z_{2}}(\log z_{2})^{k}}&\nn
&{\displaystyle =w_{(2)}\boxtimes_{P(z_{2})}\pi_{n}(w_{(1)}
\boxtimes_{P(z_{1})}w_{(3)})}.&
\end{eqnarray}
But by the definition of 
$1_{W_{2}}\boxtimes_{P(z_{2})}\mathcal{R}_{P(z_{1})}$,
\begin{eqnarray}\label{commu-2-0-2}
\lefteqn{\overline{(1_{W_{2}}\boxtimes_{P(z_{2})}\mathcal{R}_{P(z_{1})})}
(w_{(2)}\boxtimes_{P(z_{2})}\pi_{n}(w_{(1)}\boxtimes_{P(z_{1})}w_{(3)}))}\nn
&&\quad =w_{(2)}\boxtimes_{P(z_{2})}
\mathcal{R}_{P(z_{1})}
(\pi_{n}(w_{(1)}\boxtimes_{P(z_{1})}w_{(3)}))\nn
&&\quad =w_{(2)}\boxtimes_{P(z_{2})}
\pi_{n}(\mathcal{R}_{P(z_{1})}
(w_{(1)}\boxtimes_{P(z_{1})}w_{(3)}))\quad\quad
\end{eqnarray}
for $n\in \R$. {}From (\ref{commu-2-0-1}), (\ref{commu-2-0-2}) and 
the definitions of $\Y_{1}$ and $\tilde{\Y}_{1}$,
we obtain
\begin{eqnarray}\label{commu-2-0-3}
\lefteqn{\overline{(1_{W_{2}}\boxtimes_{P(z_{2})}\mathcal{R}_{P(z_{1})})}
((w_{(2)})^{\Y_{1}}_{m;k}\pi_{n}(w_{(1)}\boxtimes_{P(z_{1})}
w_{(3)}))}\nn
&&=(w_{(2)})^{\tilde{\Y}_{1}}_{m;k}
\pi_{n}(\mathcal{R}_{P(z_{1})}(w_{(1)}\boxtimes_{P(z_{1})}
w_{(3)}))\nn
&&=(w_{(2)})^{\tilde{\Y}_{1}}_{m;k}
\pi_{n}(e^{z_{1}L(-1)}(w_{(3)}\boxtimes_{P(-z_{1})}
w_{(1)}))
\end{eqnarray}
for $m\in \R$ and $k=1, \dots, N$. 
Using (\ref{commu-2-0-0}) and (\ref{commu-2-0-3}),
we obtain
\begin{eqnarray}\label{commu-2-1}
\lefteqn{\overline{(1_{W_{2}}\boxtimes_{P(z_{2})}\mathcal{R}_{P(z_{1})})}
(\Y_{1}(w_{(2)}, 
z_{2}+\epsilon)(w_{(1)}\boxtimes_{P(z_{1})}w_{(3)}))}\nn
&&=\overline{(1_{W_{2}}\boxtimes_{P(z_{2})}\mathcal{R}_{P(z_{1})})}
\nn
&&\quad\quad\left(\sum_{m, n\in \R}\sum_{k=1}^{N}
(w_{(2)})^{\Y_{1}}_{m;k}\pi_{n}(w_{(1)}\boxtimes_{P(z_{1})}
w_{(3)})e^{(-m-1)\log (z_{2}+\epsilon)}(\log (z_{2}+\epsilon))^{k}\right)\nn
&&=\sum_{m, n\in \R}\sum_{k=1}^{N}
\overline{(1_{W_{2}}\boxtimes_{P(z_{2})}\mathcal{R}_{P(z_{1})})}
((w_{(2)})^{\Y_{1}}_{m;k}\pi_{n}(w_{(1)}\boxtimes_{P(z_{1})}
w_{(3)}))\cdot\nn
&&\quad\quad\quad\quad\quad\quad\cdot
e^{(-m-1)\log (z_{2}+\epsilon)}(\log (z_{2}+\epsilon))^{k}\nn
&&=\sum_{m, n\in \R}\sum_{k=1}^{N}
(w_{(2)})^{\tilde{\Y}_{1}}_{m;k}
\pi_{n}(e^{z_{1}L(-1)}(w_{(3)}\boxtimes_{P(-z_{1})}
w_{(1)}))
e^{(-m-1)\log (z_{2}+\epsilon)}(\log (z_{2}+\epsilon))^{k}\nn
&&=\tilde{\Y}_{1}(w_{(2)}, z_{2}+\epsilon)(e^{z_{1}L(-1)}
(w_{(3)}\boxtimes_{P(-z_{1})}
w_{(1)})).
\end{eqnarray}

Let $\Y_{3}=\Y_{\boxtimes_{P(z_{1})}, 0}$ be the intertwining operator
of type
\[
{W_{2}\boxtimes_{P(z_{2}-z_{1})}(W_{1}\boxtimes_{P(-z_{1})}W_{3})\choose 
W_{2}\;\;W_{1}\boxtimes_{(-z_{1})}W_{3}}.
\]
Then by the definition of the parallel transport isomorphism,
for $m\in \R$
we have 
\begin{eqnarray}\label{commu-2-2}
\lefteqn{\overline{\mathcal{T}_{\gamma}}(
\tilde{\Y}_{1}(w_{(2)}, z_{2}+\epsilon)\pi_{m}(e^{z_{1}L(-1)}(w_{(3)}
\boxtimes_{P(-z_{1})}w_{(1)})))}\nn
&&=\overline{\mathcal{T}_{\gamma}}(
e^{\epsilon L(-1)}\tilde{\Y}_{1}(w_{(2)}, z_{2})e^{\epsilon L(-1)}
\pi_{m}(e^{z_{1}L(-1)}(w_{(3)}
\boxtimes_{P(-z_{1})}w_{(1)})))\nn
&&=\overline{\mathcal{T}_{\gamma}}(
e^{\epsilon L(-1)}w_{(2)}\boxtimes_{P(z_{2})}e^{\epsilon L(-1)}
\pi_{m}(e^{z_{1}L(-1)}(w_{(3)}
\boxtimes_{P(-z_{1})}w_{(1)})))\nn
&&=e^{\epsilon L(-1)}\Y_{3}(w_{(2)}, z_{2})e^{\epsilon L(-1)}
\pi_{m}(e^{z_{1} L(-1)}(w_{(3)}
\boxtimes_{P(-z_{1})}w_{(1)}))\nn
&&=\Y_{3}(w_{(2)}, z_{2}+\epsilon)
\pi_{m}(e^{z_{1} L(-1)}(w_{(3)}
\boxtimes_{P(-z_{1})}w_{(1)})).
\end{eqnarray}
Since $|z_{2}+\epsilon|>|z_{1}|>0$, the sums of both
sides of (\ref{commu-2-2}) over $m\in \R$ are absolutely 
convergent. Note also that $|z_{2}-z_{1}+\epsilon|>|z_{1}|>0$.
Thus we have
\begin{eqnarray}\label{commu-2-3}
\lefteqn{\sum_{m\in \R}\overline{\mathcal{T}_{\gamma}}(
\tilde{\Y}_{1}(w_{(2)}, z_{2}+\epsilon)\pi_{m}(e^{z_{1}L(-1)}(w_{(3)}
\boxtimes_{P(-z_{1})}w_{(1)})))}\nn
&&=\overline{\mathcal{T}_{\gamma}}(
\tilde{\Y}_{1}(w_{(2)}, z_{2}+\epsilon)e^{z_{1}L(-1)}(w_{(3)}
\boxtimes_{P(-z_{1})}w_{(1)}))\nn
&&=\overline{\mathcal{T}_{\gamma}}(e^{z_{1}L(-1)}
\tilde{\Y}_{1}(w_{(2)}, z_{2}-z_{1}+\epsilon)(w_{(3)}
\boxtimes_{P(-z_{1})}w_{(1)}))
\end{eqnarray}
and
\begin{eqnarray}\label{commu-2-4}
\lefteqn{\sum_{m\in \R}\Y_{3}(w_{(2)}, z_{2}+\epsilon)
\pi_{m}(e^{z_{1} L(-1)}(w_{(3)}
\boxtimes_{P(-z_{1})}w_{(1)}))}\nn
&&=\Y_{3}(w_{(2)}, z_{2}+\epsilon)
e^{z_{1} L(-1)}(w_{(3)}
\boxtimes_{P(-z_{1})}w_{(1)})\nn
&&=e^{z_{1} L(-1)}\Y_{3}(w_{(2)}, z_{2}-z_{1}+\epsilon)
(w_{(3)}
\boxtimes_{P(-z_{1})}w_{(1)}).
\end{eqnarray}
{}From (\ref{commu-2-2})--(\ref{commu-2-4}), we obtain
\begin{eqnarray}\label{commu-2-5}
\lefteqn{\overline{\mathcal{T}_{\gamma}}(e^{z_{1}L(-1)}
\Y_{1}(w_{(2)}, z_{2}-z_{1}+\epsilon)(w_{(3)}
\boxtimes_{P(-z_{1})}w_{(1)}))}\nn
&&=e^{z_{1} L(-1)}\Y_{3}(w_{(2)}, z_{2}-z_{1}+\epsilon)
(w_{(3)}
\boxtimes_{P(-z_{1})}w_{(1)}).
\end{eqnarray}
{}From (\ref{commu-2-1}) and (\ref{commu-2-5}), we obtain
\begin{eqnarray}\label{commu-2-6}
\lefteqn{\overline{\mathcal{T}_{\gamma}\circ 
(1_{W_{2}}\boxtimes_{P(z_{2})}\mathcal{R}_{P(z_{1})})}
(\Y_{1}(w_{(2)}, z_{2}+\epsilon)(w_{(1)}
\boxtimes_{P(z_{1})}w_{(3)}))}\nn
&&=e^{z_{1} L(-1)}\Y_{3}(w_{(2)}, z_{2}-z_{1}+\epsilon)
(w_{(3)}
\boxtimes_{P(-z_{1})}w_{(1)}).
\end{eqnarray}

For any 
\[
w'\in (W_{2}\boxtimes_{P(z_{2}-z_{1})}(W_{1}\boxtimes_{P(-z_{1})}W_{3}))',
\]
we have
\begin{eqnarray*}
\lefteqn{\langle w', \overline{\mathcal{T}_{\gamma}\circ 
(1_{W_{2}}\boxtimes_{P(z_{2})}\mathcal{R}_{P(z_{1})})}
(\Y_{1}(w_{(2)}, z_{2}+\epsilon)(w_{(1)}
\boxtimes_{P(z_{1})}w_{(3)}))\rangle}\nn
&&=\langle w', e^{z_{1} L(-1)}\Y_{3}(w_{(2)}, z_{2}-z_{1}+\epsilon)
(w_{(3)}
\boxtimes_{P(-z_{1})}w_{(1)})\rangle,
\end{eqnarray*}
or equivalently,
\begin{eqnarray}\label{commu-2-7}
\lefteqn{\langle ((1_{W_{2}}\boxtimes_{P(z_{2})}\mathcal{R}_{P(z_{1})})
\circ \mathcal{T}_{\gamma})(w)', \Y_{1}(w_{(2)}, z_{2}+\epsilon)(w_{(1)}
\boxtimes_{P(z_{1})}w_{(3)})\rangle}\nn
&&=\langle w', e^{z_{1} L(-1)}\Y_{3}(w_{(2)}, z_{2}-z_{1}+\epsilon)
(w_{(3)}
\boxtimes_{P(-z_{1})}w_{(1)})\rangle,
\end{eqnarray}
The left- and right-hand sides of (\ref{commu-2-7})
are values at $(\zeta_{1}, \zeta_{2})
=(z_{1}, z_{2}+\epsilon)$  of  single-valued
analytic functions of $\zeta_{1}$ and $\zeta_{2}$
defined in the region 
\[
\{(\zeta_{1}, \zeta_{2})\in \C^{2}\;|\;\zeta_{1}\ne 0,\;
\zeta_{2}\ne 0,\;\zeta_{1}\ne \zeta_{2},\;
0\le \arg \zeta_{1}, \arg \zeta_{2}, 
\arg (\zeta_{1}-\zeta_{2})<2\pi\}.
\]
Also, by the definition of tensor product of three elements
above, the values of these analytic functions at $(\zeta_{1}, \zeta_{2})
=(z_{1}, z_{2})$ are equal to 
\[
\langle ((1_{W_{2}}\boxtimes_{P(z_{2})}\mathcal{R}_{P(z_{1})})
\circ \mathcal{T}_{\gamma})(w)', w_{(2)}\boxtimes_{P(z_{2})}(w_{(1)}
\boxtimes_{P(z_{1})}w_{(3)})\rangle
\]
and 
\[
\langle w', e^{z_{1} L(-1)}w_{(2)}\boxtimes_{P(z_{2}-z_{1})}
(w_{(3)}
\boxtimes_{P(-z_{1})}w_{(1)})\rangle,
\]
respectively.
Thus we can take the limit $\epsilon\to 0$ on both sides of 
(\ref{commu-2-7}) and obtain
\begin{eqnarray}\label{commu-2-8}
\lefteqn{\langle ((1_{W_{2}}\boxtimes_{P(z_{2})}\mathcal{R}_{P(z_{1})})
\circ \mathcal{T}_{\gamma})(w)', w_{(2)}\boxtimes_{P(z_{2})}(w_{(1)}
\boxtimes_{P(z_{1})}w_{(3)})\rangle}\nn
&&\quad\quad\quad=\langle w', e^{z_{1} L(-1)}w_{(2)}\boxtimes_{P(z_{2}-z_{1})}
(w_{(3)}\boxtimes_{P(-z_{1})}w_{(1)})\rangle.
\end{eqnarray}
Since $w'$ is arbitrary, (\ref{commu-2-8}) is equivalent to 
(\ref{commu-2}).
\epfv

\subsection{The coherence properties}

We shall use the following terminology concerning tensor categories;
cf. \cite{Ma}, \cite{T} and \cite{BK}:

A {\it preadditive category} (or {\it Ab-category}) is a category in
which each hom-set is an additive abelian group such that the
composition of morphisms is bilinear.  An {\it additive category} is a
preadditive category which has a zero object and a biproduct for each
pair of objects. An {\it abelian category} is an additive category
such that every morphism has a kernel and a cokernel, every monic
morphism is a kernel and every epic morphism is a cokernel.

A {\it monoidal category} is a category $\mathcal{C}$ equipped with a
{\it monoidal} or {\it tensor product bifunctor} $\boxtimes:
\mathcal{C}\times \mathcal{C}\to \mathcal{C}$, a {\it unit object}
$V$, a {\it natural associativity isomorphism} $\mathcal{A}: \boxtimes
\circ (1_{\mathcal{C}}\times \boxtimes) \to \boxtimes \circ (\boxtimes
\times 1_{\mathcal{C}})$, a {\it natural left identity isomorphism}
$l: V\boxtimes \cdot \to 1_{\mathcal{C}}$, and a {\it natural right
identity isomorphism} $r: \cdot \boxtimes V \to 1_{\mathcal{C}}$, such
that the {\it pentagon diagram}

\begin{equation}\label{pentagon-diag}
\mbox{\rm 
\begin{picture}(150,100)(30,20)
\put(-65,20){$((W_{1}\boxtimes
W_{2})\boxtimes W_{3})
\boxtimes W_{4}$}
\put(148,20){$(W_{1}\boxtimes
(W_{2}\boxtimes W_{3}))
\boxtimes W_{4}$}
\put(-65,68){$(W_{1}\boxtimes
W_{2})\boxtimes (W_{3}
\boxtimes W_{4})$}
\put(148,68){$W_{1}\boxtimes
((W_{2}\boxtimes W_{3})
\boxtimes W_{4})$}
\put(45,113){$W_{1}\boxtimes
(W_{2}\boxtimes (W_{3}
\boxtimes W_{4}))$}

\put(145,23){\vector(-1,0){78}}
\put(0,60){\vector(0,-1){28}}
\put(210,60){\vector(0,-1){28}}

\put(80,105){\vector(-3,-1){75}}
\put(130,105){\vector(3,-1){75}}
\end{picture}}
\end{equation}
and the {\it triangle diagram}
\begin{equation}\label{trianle-diag}
\begin{CD}
(W_{1}\boxtimes V)\boxtimes W_{2}&
@>>>
&W_{1}\boxtimes (V\boxtimes W_{2})\\
@VVV&
&@VVV\\
W_{1}\boxtimes W_{2}&@>>=>
&W_{1}\boxtimes W_{2}
\end{CD}
\end{equation}
commute and such that the morphisms $l_{V}: V\boxtimes V\to V$ and
$r_{V}: V\boxtimes V\to V$ are equal.

A {\it braided monoidal category} is a monoidal category with a {\it
natural braiding isomorphism}
\[\mathcal{R}: \boxtimes \to \boxtimes \circ \sigma_{12},\]
where $\sigma_{12}$ is the permutation functor on
$\mathcal{C}\times \mathcal{C}$, such that the two hexagon diagrams
(with $\mathcal{R}^{\pm 1}$)

\vspace{3em}

\begin{picture}(200,175)(-100,0)
\put(73,165){\footnotesize $(W_{1}\boxtimes W_{2})
\boxtimes W_{3}$}
\put(-47,66){\footnotesize $(W_{2}\boxtimes W_{1})
\boxtimes W_{3}$}
\put(193,66){\footnotesize $W_{1}\boxtimes (W_{2}
\boxtimes W_{3})$}
\put(-47,-30){\footnotesize $W_{2}\boxtimes (W_{1}
\boxtimes W_{3})$}
\put(193,-30){\footnotesize $(W_{2}\boxtimes W_{3})
\boxtimes W_{1}$}
\put(73,-129){\footnotesize $W_{2}\boxtimes (W_{3}
\boxtimes W_{1}$)}

\put(-26,114){\footnotesize
$\mathcal{R}^{\pm 1}\boxtimes 1_{W_{3}}$}

\put(200,114){\footnotesize
$\A^{-1}$}

\put(-36,18){\footnotesize
$\A^{-1}$}

\put(200,-78){\footnotesize
$\A^{-1}$}

\put(-26,-78){\footnotesize
$1_{W_{2}} \boxtimes \mathcal{R}^{\pm 1}$}

\put(235,18){\footnotesize
$\mathcal{R}^{\pm 1}$}


\put(-10,60){\vector(0,-1){78}}
\put(230,60){\vector(0,-1){78}}


\put(-7,-36){\vector(1,-1){78}}
\put(225,-36){\vector(-1,-1){78}}

\put(75,155){\vector(-1,-1){78}}
\put(145,155){\vector(1,-1){78}}


\end{picture}
\vskip 1.4in
\begin{equation}\label{hexagon-diag}
\end{equation}
commute.

An {\it additive braided monoidal category} is an additive category
with a compatible braided monoidal category structure. A {\it tensor
category} is an abelian category with a compatible monoidal category
structure. A {\it braided tensor category} is an abelian category with
a compatible braided monoidal category structure, that is, a tensor
category with a compatible braiding structure.

\begin{theo}\label{main}
Let $V$ be a M\"{o}bius or conformal vertex algebra and $\mathcal{C}$
a full subcategory of $\mathcal{M}_{sg}$ or $\mathcal{GM}_{sg}$
(recall Notation \ref{MGM}) satisfying Assumptions \ref{assum-assoc},
\ref{assum-V} and \ref{assum-con}.  Then the category $\mathcal{C}$,
equipped with the tensor product bifunctor $\boxtimes$, the unit
object $V$, the braiding isomorphisms $\mathcal{R}$, the associativity
isomorphisms $\mathcal{A}$, and the left and right unit isomorphisms
$l$ and $r$, is an additive braided monoidal category.
\end{theo}
\pf We need only prove the coherence properties. We prove the
commutativity of the pentagon diagram first.  Let $W_{1}$, $W_{2}$,
$W_{3}$ and $W_{4}$ be objects of $\mathcal{C}$ and let $z_{1}, z_{2},
z_{3}\in \R$ such that
\begin{eqnarray}
&|z_{1}|>|z_{2}|>|z_{3}|>
|z_{13}|>|z_{23}|>|z_{12}|>0,&\nn
&|z_{1}|>|z_{3}|+|z_{23}|>0,&\nn
&|z_{2}|>|z_{12}|+|z_{3}|>0,&\nn
&|z_{3}|>|z_{23}|+|z_{12}|>0,&
\end{eqnarray}
where $z_{12}=z_{1}-z_{2}$, $z_{13}=z_{1}-z_{3}$ and $z_{23}=z_{2}-z_{3}$.
For example, we can take $z_{1}=7$, $z_{2}=6$ and $z_{3}=4$.
We first prove the commutativity of the following diagram:

\vspace{3.5em}

\begin{picture}(150,100)(-70,0)

\put(-85,20){\footnotesize $((W_{1}\boxtimes_{P(z_{12})}
W_{2})\boxtimes_{P(z_{23})} W_{3})
\boxtimes_{P(z_{3})} W_{4}$}
\put(123,20){\footnotesize $(W_{1}\boxtimes_{P(z_{13})}
(W_{2}\boxtimes_{P(z_{23})} W_{3}))
\boxtimes_{P(z_{3})} W_{4}$}
\put(-85,68){\footnotesize $(W_{1}\boxtimes_{P(z_{12})}
W_{2})\boxtimes_{P(z_{2})} (W_{3}
\boxtimes_{P(z_{3})} W_{4})$}
\put(123,68){\footnotesize $W_{1}\boxtimes_{P(z_{1})}
((W_{2}\boxtimes_{P(z_{23})} W_{3})
\boxtimes_{P(z_{3})} W_{4})$}
\put(20,116){\footnotesize $W_{1}\boxtimes_{P(z_{1})}
(W_{2}\boxtimes_{P(z_{2})} (W_{3}
\boxtimes_{P(z_{3})} W_{4})).$}

\put(120,23){\vector(-1,0){25}}
\put(5,60){\vector(0,-1){28}}
\put(212,60){\vector(0,-1){28}}

\put(100,105){\vector(-3,-1){75}}
\put(110,105){\vector(3,-1){75}}
\end{picture}
\begin{equation}\label{pent1}
\end{equation}
For $w_{(1)}\in W_{1}$, $w_{(2)}\in W_{2}$, $w_{(3)}\in W_{3}$ and
$w_{(4)}\in W_{4}$, we consider
\[
w_{(1)}\boxtimes_{P(z_{1})} (w_{(2)}\boxtimes_{P(z_{2})}
(w_{(3)}\boxtimes_{P(z_{3})}w_{(4)}))\in
\overline{W_{1}\boxtimes_{P(z_{1})} (W_{2}\boxtimes_{P(z_{2})}
(W_{3}\boxtimes_{P(z_{3})}W_{4}))}.
\]
By the characterizations of the associativity isomorphisms,
we see that the compositions of the natural extensions
of the module maps in the two routes in (\ref{pent1}) applied to
this element both give
\[
((w_{(1)}\boxtimes_{P(z_{12})} w_{(2)})\boxtimes_{P(z_{23})}
w_{(3)})\boxtimes_{P(z_{3})}w_{(4)}\in
\overline{((W_{1}\boxtimes_{P(z_{12})} W_{2})\boxtimes_{P(z_{23})}
W_{3})\boxtimes_{P(z_{3})}W_{4}}.
\]
Since the homogeneous components of
\[
w_{(1)}\boxtimes_{P(z_{1})} (w_{(2)}\boxtimes_{P(z_{2})}
(w_{(3)}\boxtimes_{P(z_{3})}w_{(4)}))
\]
for $w_{(1)}\in W_{1}$, $w_{(2)}\in W_{2}$, $w_{(3)}\in W_{3}$ and
$w_{(4)}\in W_{4}$ span
\[
W_{1}\boxtimes_{P(z_{1})} (W_{2}\boxtimes_{P(z_{2})}
(W_{3}\boxtimes_{P(z_{3})}W_{4})),
\]
the diagram (\ref{pent1}) above commutes.

On the other hand, by the definition of $\A$, the  diagrams
\begin{equation}\label{pent2}
\begin{picture}(60,90)(20,0)
\put(-145,68){\footnotesize $W_{1}\boxtimes_{P(z_{1})}
(W_{2}\boxtimes_{P(z_{2})} (W_{3}
\boxtimes_{P(z_{3})} W_{4}))$}
\put(63,68){\footnotesize $(W_{1}\boxtimes_{P(z_{12})}
W_{2})\boxtimes_{P(z_{2})} (W_{3}
\boxtimes_{P(z_{3})} W_{4})$}
\put(-105,20){\footnotesize $W_{1}\boxtimes
(W_{2}\boxtimes (W_{3}
\boxtimes W_{4}))$}
\put(88,20){\footnotesize $(W_{1}\boxtimes
W_{2})\boxtimes (W_{3}
\boxtimes W_{4})$}

\put(5,23){\vector(1,0){78}}
\put(30,71){\vector(1,0){25}}
\put(-55,60){\vector(0,-1){28}}
\put(152,60){\vector(0,-1){28}}
\end{picture}
\end{equation}
\begin{equation}\label{pent3}
\begin{picture}(60,90)(20,0)
\put(-145,68){\footnotesize $(W_{1}\boxtimes_{P(z_{12})}
W_{2})\boxtimes_{P(z_{2})} (W_{3}
\boxtimes_{P(z_{3})} W_{4})$}
\put(63,68){\footnotesize $((W_{1}\boxtimes_{P(z_{12})}
W_{2})\boxtimes_{P(z_{23})} W_{3})
\boxtimes_{P(z_{3})} W_{4}$}
\put(-105,20){\footnotesize $(W_{1}\boxtimes
W_{2})\boxtimes (W_{3}
\boxtimes W_{4})$}
\put(88,20){\footnotesize $((W_{1}\boxtimes
W_{2})\boxtimes W_{3})
\boxtimes W_{4}$}

\put(5,23){\vector(1,0){78}}
\put(34,71){\vector(1,0){25}}
\put(-55,60){\vector(0,-1){28}}
\put(152,60){\vector(0,-1){28}}
\end{picture}
\end{equation}
\begin{equation}\label{pent4}
\begin{picture}(60,90)(20,0)
\put(-145,68){\footnotesize $W_{1}\boxtimes_{P(z_{1})}
(W_{2}\boxtimes_{P(z_{2})} (W_{3}
\boxtimes_{P(z_{3})} W_{4}))$}
\put(63,68){\footnotesize $W_{1}\boxtimes_{P(z_{1})}
((W_{2}\boxtimes_{P(z_{23})} W_{3})
\boxtimes_{P(z_{3})} W_{4})$}
\put(-105,20){\footnotesize $W_{1}\boxtimes
(W_{2}\boxtimes (W_{3}
\boxtimes W_{4}$))}
\put(88,20){\footnotesize $W_{1}\boxtimes
((W_{2}\boxtimes W_{3})
\boxtimes W_{4})$}

\put(5,23){\vector(1,0){78}}
\put(34,71){\vector(1,0){25}}
\put(-55,60){\vector(0,-1){28}}
\put(152,60){\vector(0,-1){28}}
\end{picture}
\end{equation}
\begin{equation}\label{pent5}
\begin{picture}(60,80)(20,0)
\put(-145,68){\footnotesize $W_{1}\boxtimes_{P(z_{1})}
((W_{2}\boxtimes_{P(z_{23})} W_{3})
\boxtimes_{P(z_{3})} W_{4}))$}
\put(63,68){\footnotesize $(W_{1}\boxtimes_{P(z_{13})}
(W_{2}\boxtimes_{P(z_{23})} W_{3}))
\boxtimes_{P(z_{3})} W_{4}$}

\put(-105,20){\footnotesize $W_{1}\boxtimes
((W_{2}\boxtimes W_{3})
\boxtimes W_{4})$}
\put(88,20){\footnotesize $(W_{1}\boxtimes
(W_{2}\boxtimes W_{3}))
\boxtimes W_{4}$}

\put(5,23){\vector(1,0){78}}
\put(34,71){\vector(1,0){25}}
\put(-55,60){\vector(0,-1){28}}
\put(152,60){\vector(0,-1){28}}
\end{picture}
\end{equation}
\begin{equation}\label{pent6}
\begin{picture}(60,90)(20,0)
\put(-145,68){\footnotesize $(W_{1}\boxtimes_{P(z_{13})}
(W_{2}\boxtimes_{P(z_{23})} W_{3}))
\boxtimes_{P(z_{3})} W_{4}$}
\put(63,68){\footnotesize $((W_{1}\boxtimes_{P(z_{12})}
W_{2})\boxtimes_{P(z_{23})} W_{3})
\boxtimes_{P(z_{3})} W_{4}$}
\put(-105,20){\footnotesize $(W_{1}\boxtimes
(W_{2}\boxtimes W_{3}))
\boxtimes W_{4}$}
\put(88,20){\footnotesize $((W_{1}\boxtimes
W_{2})\boxtimes W_{3})
\boxtimes W_{4}$}

\put(5,23){\vector(1,0){78}}
\put(34,71){\vector(1,0){25}}
\put(-55,60){\vector(0,-1){28}}
\put(152,60){\vector(0,-1){28}}
\end{picture}
\end{equation}
all commute.  Combining all the diagrams (\ref{pent1})--(\ref{pent6})
above, we see that the pentagon diagram

\vspace{3.5em}

\begin{picture}(150,100)(-70,0)

\put(-45,20){\footnotesize $((W_{1}\boxtimes
W_{2})\boxtimes W_{3})
\boxtimes W_{4}$}
\put(148,20){\footnotesize $(W_{1}\boxtimes
(W_{2}\boxtimes W_{3}))
\boxtimes W_{4}$}
\put(-45,68){\footnotesize $(W_{1}\boxtimes
W_{2})\boxtimes (W_{3}
\boxtimes W_{4})$}
\put(148,68){\footnotesize $W_{1}\boxtimes
((W_{2}\boxtimes W_{3})
\boxtimes W_{4})$}
\put(55,113){\footnotesize $W_{1}\boxtimes
(W_{2}\boxtimes (W_{3}
\boxtimes W_{4}))$}

\put(145,23){\vector(-1,0){78}}
\put(10,60){\vector(0,-1){28}}
\put(202,60){\vector(0,-1){28}}

\put(100,105){\vector(-3,-1){75}}
\put(110,105){\vector(3,-1){75}}
\end{picture}

\noindent also commutes.

Next we prove the commutativity of the hexagon diagrams.  We prove
only the commutativity of the hexagon diagram involving $\mathcal{R}$;
the proof of the commutativity of the other hexagon diagram is the
same.  Let $W_{1}$, $W_{2}$ and $W_{3}$ be objects of $\mathcal{C}$
and let $z_{1}, z_{2}\in \C^{\times}$ satisfying
$|z_{1}|=|z_{2}|=|z_{1}-z_{2}|$ and let $z_{12}=z_{1}-z_{2}$.  We
first prove the commutativity of the following diagram:


\begin{picture}(200,175)(-100,0)
\put(50,162){\footnotesize $(W_{1}\boxtimes_{P(z_{12})} W_{2})
\boxtimes_{P(z_{2})} W_{3}$}
\put(-70,66){\footnotesize $(W_{2}\boxtimes_{P(-z_{12})} W_{1})
\boxtimes_{P(z_{2})} W_{3}$}
\put(170,66){\footnotesize $W_{1}\boxtimes_{P(z_{1})} (W_{2}
\boxtimes_{P(z_{2})} W_{3})$}
\put(-70,18){\footnotesize $(W_{2}\boxtimes_{P(-z_{12})} W_{1})
\boxtimes_{P(z_{1})} W_{3}$}
\put(-70,-30){\footnotesize $W_{2}\boxtimes_{P(z_{2})} (W_{1}
\boxtimes_{P(z_{1})} W_{3})$}
\put(170,-30){\footnotesize $(W_{2}\boxtimes_{P(z_{2})} W_{3})
\boxtimes_{P(-z_{1})} W_{1}$}
\put(120,-78){\footnotesize $
W_{2}\boxtimes_{P(-z_{12})} (W_{3}
\boxtimes_{P(-z_{1})} W_{1})$}
\put(50,-126){\footnotesize $W_{2}\boxtimes_{P(z_{2})} (W_{3}
\boxtimes_{P(-z_{1})} W_{1}$)}

\put(-70,114){\footnotesize
$\mathcal{R}_{P(z_{12})}\boxtimes_{P(z_{2})} 1_{W_{3}}$}

\put(200,114){\footnotesize
$\left(\A_{P(z_{1}), P(z_{2})}
^{P(z_{12}), P(z_{2})}\right)^{-1}$}

\put(-100,-5){\footnotesize
$\left(\A_{P(z_{2}), P(z_{1})}
^{P(-z_{12}), P(z_{1})}\right)^{-1}$}

\put(220,-55){\footnotesize
$\left(\A_{P(-z_{12}), P(-z_{1})}
^{P(z_{2}), P(-z_{1})}\right)^{-1}$}

\put(-60,-78){\footnotesize
$1_{W_{2}} \boxtimes_{P(z_{2})} \mathcal{R}_{P(z_{1})}$}

\put(235,18){\footnotesize
$\mathcal{R}_{P(z_{1})}$}

\put(-35,43){\footnotesize
$\mathcal{T}_{\gamma_{1}}$}

\put(170,-105){\footnotesize
$\mathcal{T}_{\gamma_{2}}$}

\put(-10,60){\vector(0,-1){28}}
\put(230,60){\vector(0,-1){78}}

\put(-10,12){\vector(0,-1){28}}

\put(-7,-36){\vector(1,-1){78}}
\put(225,-36){\vector(-1,-1){30}}

\put(75,155){\vector(-1,-1){78}}
\put(145,155){\vector(1,-1){78}}

\put(175,-85){\vector(-1,-1){30}}

\end{picture}
\vskip 1.4in
\begin{equation}\label{hexagon1}
\end{equation}
where $\gamma_{1}$ and $\gamma_{2}$ are paths {}from $z_{2}$ to $z_{1}$
and {}from $-z_{12}$ to $z_{2}$, respectively, in $\C$ with a cut along
the nonnegative real line.

Let $w_{(1)}\in W_{1}$, $w_{(2)}\in W_{2}$ and $w_{(3)}\in W_{3}$.  By
the results proved in the preceding section, we see that the images
of the element
\[
(w_{(1)}\boxtimes_{P(z_{12})}
w_{(2)})\boxtimes_{P(z_{2})}w_{(3)}
\]
under the natural extension to
\[
\overline{(W_{1}\boxtimes_{P(z_{12})} W_{2})
\boxtimes_{P(z_{2})} W_{3}}
\]
of the compositions of the maps in both the left and right routes in
(\ref{hexagon1}) {}from 
\[
(W_{1}\boxtimes_{P(z_{12})} W_{2})
\boxtimes_{P(z_{2})} W_{3}
\]
to 
\[
W_{2}\boxtimes_{P(z_{2})} (W_{3}
\boxtimes_{P(-z_{1})} W_{1})
\]
are 
\[
w_{(2)}\boxtimes_{P(z_{2})}
(e^{z_{1}L(-1)}(w_{(3)}\boxtimes_{P(-z_{1})}w_{(1)})).
\]
Since the homogeneous components of 
\[
(w_{(1)}\boxtimes_{P(z_{12})}
w_{(2)})\boxtimes_{P(z_{2})}w_{(3)}
\]
for $w_{(1)}\in W_{1}$, $w_{(2)}\in W_{2}$ and $w_{(3)}\in W_{3}$
span
\[
(W_{1}\boxtimes_{P(z_{12})} W_{2})
\boxtimes_{P(z_{2})} W_{3},
\]
the diagram (\ref{hexagon1}) commutes.

Now we consider the following diagrams:
\begin{equation}\label{hexagon2}
\begin{CD}
(W_{1}\boxtimes_{P(z_{12})} W_{2})
\boxtimes_{P(z_{2})} W_{3}&@>>>&(W_{1}\boxtimes W_{2})
\boxtimes W_{3}\\
@VVV&&@VVV\\
(W_{2}\boxtimes_{P(-z_{12})} W_{1})
\boxtimes_{P(z_{2})} W_{3}&@>>>&(W_{2}\boxtimes W_{1})
\boxtimes W_{3}
\end{CD}
\end{equation}
\vspace{1.5em}
\begin{equation}\label{hexagon3}
\begin{CD}
(W_{2}\boxtimes_{P(-z_{12})} W_{1})
\boxtimes_{P(z_{2})} W_{3}&@>>>&(W_{2}\boxtimes W_{1})
\boxtimes W_{3}\\
@VVV&&@VVV\\
(W_{2}\boxtimes_{P(-z_{12})} W_{1})
\boxtimes_{P(z_{1})} W_{3}&@>>>&(W_{2}\boxtimes W_{1})
\boxtimes W_{3}
\end{CD}
\end{equation}
\vspace{1.5em}
\begin{equation}\label{hexagon4}
\begin{CD}
(W_{2}\boxtimes_{P(-z_{12})} W_{1})
\boxtimes_{P(z_{1})} W_{3}&@>>>&(W_{2}\boxtimes W_{1})
\boxtimes W_{3}\\
@VVV&&@VVV\\
W_{2}\boxtimes_{P(z_{2})} (W_{1}
\boxtimes_{P(z_{1})} W_{3})&@>>>&W_{2}\boxtimes (W_{1}
\boxtimes W_{3})\\
\end{CD}
\end{equation}
\vspace{1.5em}
\begin{equation}\label{hexagon5}
\begin{CD}
W_{2}\boxtimes_{P(z_{2})} (W_{1}
\boxtimes_{P(z_{1})} W_{3})&@>>>&W_{2}\boxtimes (W_{1}
\boxtimes W_{3})\\
@VVV&&@VVV\\
W_{2}\boxtimes_{P(z_{2})} (W_{3}
\boxtimes_{P(-z_{1})} W_{1})&@>>>&W_{2}\boxtimes (W_{3}
\boxtimes W_{1})
\end{CD}
\end{equation}
\vspace{1.5em}
\begin{equation}\label{hexagon6}
\begin{CD}
(W_{1}\boxtimes_{P(z_{12})} W_{2})
\boxtimes_{P(z_{2})} W_{3}&@>>>&(W_{1}\boxtimes W_{2})
\boxtimes W_{3}\\
@VVV&&@VVV\\
W_{1}\boxtimes_{P(z_{1})} (W_{2}
\boxtimes_{P(z_{2})} W_{3})&@>>>&W_{1}\boxtimes (W_{1}
\boxtimes W_{3})
\end{CD}
\end{equation}
\vspace{1.5em}
\begin{equation}\label{hexagon7}
\begin{CD}
W_{1}\boxtimes_{P(z_{1})} (W_{2}
\boxtimes_{P(z_{2})} W_{3})&@>>>&W_{1}\boxtimes (W_{1}
\boxtimes W_{3})\\
@VVV&&@VVV\\
(W_{2}\boxtimes_{P(z_{2})} W_{3})
\boxtimes_{P(-z_{1})} W_{1}&@>>>&(W_{2}\boxtimes W_{3})
\boxtimes W_{1}
\end{CD}
\end{equation}
\vspace{1.5em}
\begin{equation}\label{hexagon8}
\begin{CD}
(W_{2}\boxtimes_{P(z_{2})} W_{3})
\boxtimes_{P(-z_{1})} W_{1}&@>>>&(W_{2}\boxtimes W_{3})
\boxtimes W_{1}\\
@VVV&&@VVV\\
W_{2}\boxtimes_{P(-z_{12})} (W_{3}
\boxtimes_{P(-z_{1})} W_{1})&@>>>&W_{2}\boxtimes (W_{3}
\boxtimes W_{1})
\end{CD}
\end{equation}
\vspace{1.5em}
\begin{equation}\label{hexagon9}
\begin{CD}
W_{2}\boxtimes_{P(-z_{12})} (W_{3}
\boxtimes_{P(-z_{1})} W_{1})&@>>>&W_{2}\boxtimes (W_{3}
\boxtimes W_{1})\\
@VVV&&@VVV\\
W_{2}\boxtimes_{P(z_{2})} (W_{3}
\boxtimes_{P(-z_{1})} W_{1})&@>>>&W_{2}\boxtimes (W_{3}
\boxtimes W_{1})
\end{CD}
\end{equation}
The commutativity of the diagrams (\ref{hexagon2}), (\ref{hexagon5})
and (\ref{hexagon7}) follows {}from the definition of the
commutativity isomorphism for the braided tensor category structure
and the naturality of the parallel transport isomorphisms.  The
commutativity of (\ref{hexagon4}), (\ref{hexagon6}) and
(\ref{hexagon8}) follows {}from the definition of the associativity
isomorphism for the braided tensor product structure.  The
commutativity of (\ref{hexagon3}) and (\ref{hexagon9}) follows {}from
the facts that compositions of parallel transport isomorphisms are
equal to the parallel transport isomorphisms associated to the
products of the paths and that parallel transport isomorphisms
associated to homotopically equivalent paths are equal. The
commutativity of the hexagon diagram involving $\mathcal{R}$ follows
{}from (\ref{hexagon1})--(\ref{hexagon9}).

Finally, we prove the commutativity of the triangle diagram for the
unit isomorphisms. Let $z_{1}$ and $z_{2}$ be complex numbers such
that $|z_{1}|>|z_{2}|>|z_{1}-z_{2}|>0$ and let
$z_{12}=z_{1}-z_{2}$. Also let $\gamma$ be a path {}from $z_{2}$ to
$z_{1}$ in $\C$ with a cut along the nonnegative real line. We first
prove the commutativity of the following diagram:
\begin{equation}\label{unit1}
\begin{CD}
(W_{1}\boxtimes_{P(z_{12})} V)
\boxtimes_{P(z_{2})} W_{2}&
@>(\mathcal{A}_{P(z_{1}), P(z_{2})}^{P(z_{12}), P(z_{2})})^{-1}>>
&W_{1}\boxtimes_{P(z_{1})} (V
\boxtimes_{P(z_{2})} W_{2})\\
@Vr_{z_{12};W_{1}}\boxtimes_{P(z_{2})} 1_{W_{2}}VV&
&@VV1_{W_{1}}\boxtimes_{P(z_{1})} l_{W_{2}}V\\
W_{1}\boxtimes_{P(z_{2})} W_{2}&@>>\mathcal{T}_{\gamma}>
&W_{1}\boxtimes_{P(z_{1})} W_{2}.
\end{CD}
\end{equation}

Let $w_{(1)}\in W_{1}$ and $w_{(2)}\in W_{2}$. Then we have 
\begin{eqnarray}\label{unit2}
\lefteqn{\overline{(1_{W_{1}}\boxtimes_{P(z_{1})} l_{z_{2};W_{2}})\circ
(\mathcal{A}_{P(z_{1}), P(z_{2})}^{P(z_{1}-z_{2}), P(z_{2})})^{-1}}
((w_{(1)}\boxtimes_{P(z_{12})}\mathbf{1})\boxtimes_{P(z_{2})}w_{(2)})}\nn
&&\quad\quad\quad=\overline{(1_{W_{1}}\boxtimes_{P(z_{1})} l_{z_{2};W_{2}})}
(w_{(1)}\boxtimes_{P(z_{1})}(\mathbf{1}\boxtimes_{P(z_{2})}w_{(2)})\nn
&&\quad\quad\quad=w_{(1)}\boxtimes_{P(z_{1})}w_{(2)}.
\end{eqnarray}
But 
\begin{eqnarray}\label{unit3}
\lefteqn{\overline{r_{z_{12};W_{1}}\boxtimes_{P(z_{2})} 1_{W_{2}}}
((w_{(1)}\boxtimes_{P(z_{12})}\mathbf{1})\boxtimes_{P(z_{2})}w_{(2)})}\nn
&&=(e^{z_{12}L(-1)}w_{(1)})\boxtimes_{P(z_{2})}w_{(2)}.
\end{eqnarray}
Let $\Y=\Y_{\boxtimes_{P(z_{1})}, 0}$ be the intertwining 
operator of type ${W_{1}\boxtimes_{P(z_{1})}W_{2}\choose
W_{1}\;\;W_{2}}$ corresponding to the $P(z_{1})$-intertwining 
map $\boxtimes_{P(z_{1})}$. Then by the definition of 
the parallel transport isomorphism and the $L(-1)$-derivative
property for intertwining operators, we have
\begin{eqnarray}\label{unit4}
\overline{\mathcal{T}_{\gamma}}
((e^{z_{12}L(-1)}w_{(1)})\boxtimes_{P(z_{2})}w_{(2)})
&=&\Y(e^{z_{12}L(-1)}w_{(1)}, z_{2})w_{(3)}\nn
&=&\Y(w_{(1)}, z_{1})w_{(3)}\nn
&=&w_{(1)}\boxtimes_{P(z_{1})}w_{(2)}.
\end{eqnarray}
Since the elements 
$(w_{(1)}\boxtimes_{P(z_{12})}\mathbf{1})\boxtimes_{P(z_{2})}w_{(3)}$
for $w_{(1)}\in W_{1}$ and $w_{(2)}\in W_{2}$
span $(W_{1}\boxtimes_{P(z_{12})} V)\boxtimes_{P(z_{2})}W_{3}$,
(\ref{unit2})--(\ref{unit4}) give the commutativity 
of (\ref{unit1}). 

Let $\gamma_{1}$ be a path {}from $z_{1}$ to $1$ in $\C$ with a cut along 
the nonnegative real line.  Let $\gamma_{2}$ be the product of $\gamma$ 
and $\gamma_{1}$. In particular, $\gamma_{2}$ is a path 
{}from $z_{2}$ to $1$ in $\C$ with a cut along 
the nonnegative real line. Also let $\gamma_{12}$ be a path 
{}from $z_{12}=z_{1}-z_{2}$ to $1$ in $\C$ with a cut along 
the nonnegative real line. Then we have the following commutative 
diagrams:
\begin{equation}\label{unit5}
\begin{CD}
(W_{1}\boxtimes V)\boxtimes W_{2}&
@>\mathcal{A}^{-1}>>
&W_{1}\boxtimes (V\boxtimes W_{2})\\
@V\mathcal{T}^{-1}_{\gamma_{2}}\circ (\mathcal{T}^{-1}_{\gamma_{12}}
\boxtimes_{P(z_{2})}
1_{W_{2}})VV&
&@VV\mathcal{T}^{-1}_{\gamma_{1}}\circ (1_{W_{1}}\boxtimes_{P(z_{1})}
\mathcal{T}^{-1}_{\gamma_{2}})V\\
(W_{1}\boxtimes_{P(z_{12})} V)
\boxtimes_{P(z_{2})} W_{2}&@>>
(\mathcal{A}_{P(z_{1}), P(z_{2})}^{P(z_{12}), P(z_{2})})^{-1}>
&W_{1}\boxtimes_{P(z_{1})} (V
\boxtimes_{P(z_{2})} W_{2}).
\end{CD}
\end{equation}
\begin{equation}\label{unit6}
\begin{CD}
(W_{1}\boxtimes V)
\boxtimes W_{2}&
@>\mathcal{T}_{\gamma_{2}}^{-1}\circ (\mathcal{T}_{\gamma_{12}}^{-1}
\boxtimes 1_{W_{2}})>>
&(W_{1}\boxtimes_{P(z_{12})} V)
\boxtimes_{P(z_{2})} W_{2}\\
@Vr_{W_{1}}\boxtimes 1_{W_{2}}VV&
&@VVr_{W_{1}}\boxtimes_{P(z_{2})} 1_{W_{2}}V\\
W_{1}\boxtimes W_{2}&@>>\mathcal{T}_{\gamma_{2}}^{-1}>
&W_{1}\boxtimes_{P(z_{2})} W_{2}.
\end{CD}
\end{equation}
\begin{equation}\label{unit7}
\begin{CD}
W_{1}\boxtimes_{P(z_{1})} (W_{2}
\boxtimes_{P(z_{2})} W_{2})&
@>\mathcal{T}_{\gamma_{1}}\circ (1_{W_{1}}\boxtimes_{P(z_{1})}
\mathcal{T}_{\gamma_{2}})>>
&W_{1}\boxtimes (V
\boxtimes W_{2})\\
@V1_{W_{1}}\boxtimes_{P(z_{1})} l_{W_{2}}VV&
&@VV1_{W_{1}}\boxtimes l_{W_{2}}V\\
W_{1}\boxtimes_{P(z_{1})} W_{2}&@>>\mathcal{T}_{\gamma_{1}}>
&W_{1}\boxtimes W_{2}.
\end{CD}
\end{equation}
\begin{equation}\label{unit8}
\begin{CD}
W_{1}\boxtimes_{P(z_{2})} W_{2}&
@>\mathcal{T}_{\gamma}>>
&W_{1}\boxtimes_{P(z_{1})} W_{2}\\
@V\mathcal{T}_{\gamma_{2}}VV&
&@VV\mathcal{T}_{\gamma_{1}}V\\
W_{1}\boxtimes W_{2}&&=&
&W_{1}\boxtimes W_{2}.
\end{CD}
\end{equation}
The commutativity of (\ref{unit5}) follows {}from the definition of
$\mathcal{A}$. The commutativity of (\ref{unit6}) and (\ref{unit7})
follows {}from the definition of the left and right unit isomorphisms
and the parallel transport isomorphisms. The commutativity of
(\ref{unit8}) follows {}from the fact that $\gamma_{2}$ is the product
of $\gamma$ and $\gamma_{1}$. Combining (\ref{unit1}) and
(\ref{unit5})--(\ref{unit8}), we obtain the commutativity of the
triangle diagram for the unit isomorphisms.

It is clear {}from the definition that $l_{V}=r_{V}$. 

Thus we have proved that the category $\mathcal{C}$ equipped with the
data given in Section 12.2 is a braided monoidal category.  \epfv

In the case that $\mathcal{C}$ is an abelian category, we have:

\begin{corol}
If the category $\mathcal{C}$ is an abelian category, then
$\mathcal{C}$, equipped with the tensor product bifunctor $\boxtimes$,
the unit object $V$, the braiding isomorphism $\mathcal{R}$, the
associativity isomorphism $\mathcal{A}$, and the left and right unit
isomorphisms $l$ and $r$, is a braided tensor category. \epf
\end{corol}

\begin{rema}
{\rm As we mentioned at the beginning of this section, this braided
tensor category structure has ``forgotten'' the underlying
complex-analytic vertex-tensor-categorical structure that has in fact
been developed in this work, retaining only its ``topological'' part,
but our proof has needed the vertex-tensor-categorical structure,
essentially because iterated tensor products of triples of elements
are not defined in the braided tensor category structure (recall
Remark \ref{notensorprodelts}).  Also, our category has a
contragredient functor (which we have been using extensively),
although in the present generality, we do not have rigidity.}
\end{rema}


\bigskip

\noindent {\small \sc Department of Mathematics, Rutgers University,
Piscataway, NJ 08854 (permanent address)}

\noindent {\it and}

\noindent {\small \sc Beijing International Center for Mathematical Research,
Peking University, Beijing, China}

\noindent {\em E-mail address}: yzhuang@math.rutgers.edu

\vspace{1em}

\noindent {\small \sc Department of Mathematics, Rutgers University,
Piscataway, NJ 08854}

\noindent {\em E-mail address}: lepowsky@math.rutgers.edu

\vspace{1em}

\noindent {\small \sc Department of Mathematics, Rutgers University,
Piscataway, NJ 08854}

\noindent {\em E-mail address}: linzhang@math.rutgers.edu

\end{document}